\documentclass[a4paper]{amsart} 
\usepackage[ascii]{inputenc} 

\usepackage{amssymb}
\usepackage{mathrsfs}
\usepackage[hidelinks]{hyperref}

\usepackage{xcolor}

\usepackage[shortlabels]{enumitem}
\setlist[enumerate,1]{label={(\Alph*)}}
\setlist[enumerate,2]{label={(\alph*)}}
\setlist[enumerate,3]{label={$\bullet_{\arabic*}$}}

\newenvironment{PROOF}[2][\proofname.]
{\begin{proof}[#1]\renewcommand{\qedsymbol}{\textsquare$_{\rm #2}$}}
{\end{proof}}

\newtheorem{theorem}{Theorem}[section] 
\newtheorem{claim}[theorem]{Claim}
\newtheorem{lemma}[theorem]{Lemma} 
 
\newtheorem{conclusion}[theorem]{Conclusion}

\theoremstyle{definition}
\newtheorem{convention}[theorem]{Convention}
\newtheorem{definition}[theorem]{Definition}

\newtheorem{fact}[theorem]{Fact}

\theoremstyle{remark}
\newtheorem{remark}[theorem]{Remark}

\newtheorem{notation}[theorem]{Notation}

\newcommand{\aht}{\mathrm{aht}}
\newcommand{\eht}{\mathrm{eht}}

\newcommand{\ahtn}{\mathrm{ahtn}}
\newcommand{\ehtn}{\mathrm{ehtn}}
\newcommand{\htn}{\mathrm{htn}}

\newcommand{\Cd}{\mathrm{Cd}}
\newcommand{\Col}{\mathrm{Col}}
\newcommand{\Eh}{\mathrm{Eh}}
\newcommand{\fe}{\mathrm{fe}}
\newcommand{\fs}{\mathrm{fs}}
\newcommand{\IN}{\mathrm{IN}}
\newcommand{\MA}{\mathrm{MA}}
\newcommand{\np}{\mathrm{np}}
\newcommand{\NPr}{\mathrm{NPr}}
\newcommand{\Per}{\mathrm{Per}}
\newcommand{\TrEmb}{\mathrm{TrEmb}}
\newcommand{\uq}{\mathrm{uq}}

\newcommand{\clp}{\mathrm{clp}}
\newcommand{\wsp}{\mathrm{wsp}}

\newcommand{\stg}{\mathsf{stg}}

\newcommand{\lesdot}{\mathrel{\mathord{<}\!\!\raise 0.8pt\hbox{$\scriptstyle\circ$}}}
\newcommand{\leseqdot}{\mathrel{\mathord{\leq}\!\!\raise0.8 pt\hbox{$\scriptstyle\circ$}}}

\newcommand{\hgt}{\mathrm{ht}}

\newcommand{\essp}{\mathrm{sp}}

\newcommand{\tthh}{\mathrm{th}}

\newcommand{\acc}{\mathrm{acc}}

\newcommand{\cf}{\mathrm{cf}}

\newcommand{\dom}{\mathrm{dom}}
\newcommand{\GCH}{\mathsf{GCH}}

\newcommand{\OP}{\mathrm{OP}}

\newcommand{\otp}{\mathrm{otp}}

\newcommand{\Rang}{\mathrm{Rang}}

\newcommand{\supminus}{{\text{\hspace{.05cm}--}}}

\newcommand{\SP}{\mathrm{SP}}

\newcommand{\cn}{\mathrm{cn}}

\newcommand{\bfD}{\mathbf{D}}
\newcommand{\bfd}{\mathbf{d}}

\newcommand{\bfG}{\mathbf{G}}

\newcommand{\bfh}{\mathbf{h}}
\newcommand{\bfI}{\mathbf{I}}

\newcommand{\bfK}{\mathbf{K}}

\newcommand{\bfL}{\mathbf{L}}

\newcommand{\bfn}{\mathbf{n}}

\newcommand{\bfq}{\mathbf{q}}

\newcommand{\bfT}{\mathbf{T}}
\newcommand{\bft}{\mathbf{t}}

\newcommand{\bfV}{\mathbf{V}}

\newcommand{\bbP}{\mathbb{P}}
\newcommand{\bbQ}{\mathbb{Q}}
\newcommand{\bbR}{\mathbb{R}}

\newcommand{\mn}{\medskip\noindent}
\newcommand{\sn}{\smallskip\noindent}
\newcommand{\bn}{\bigskip\noindent}

\newcommand{\clH}{\mathcal{H}}

\newcommand{\gA}{\mathfrak{A}}

\newcommand{\gD}{\mathfrak{D}}
\newcommand{\gK}{\mathfrak{K}}

\newcommand{\eps}{\varepsilon}
\newcommand{\lh}{{\ell\kern-.27ex g}}
\newcommand{\rest}{\restriction}

\newcommand{\caret}{{\char 94}}
\newcommand{\LL}{\langle}
\newcommand{\RR}{\rangle}

\newcommand{\overbar}[1]{\mkern 1.5mu\overline{\mkern-1.5mu#1\mkern-1.5mu}\mkern 1.5mu}

\newcommand{\olsi}[1]{\,\overline{\!{#1}}} 

\newcommand*{\defeq}{\mathrel{\vcenter{\baselineskip0.5ex \lineskiplimit0pt\hbox{\scriptsize.}\hbox{\scriptsize.}}}=}

\newcount\skewfactor
\def\mathunderaccent#1#2 {\let\theaccent#1\skewfactor#2
\mathpalette\putaccentunder}
\def\putaccentunder#1#2{\oalign{$#1#2$\crcr\hidewidth
\vbox to.2ex{\hbox{$#1\skew\skewfactor\theaccent{}$}\vss}\hidewidth}}
\def\name{\mathunderaccent\tilde-3 }
\def\Name{\mathunderaccent\widetilde-3 }

\newbox\noforkbox \newdimen\forklinewidth
\setbox0\hbox{$\textstyle\bigcup$}
\setbox1\hbox to \wd0{\hfil\vrule width \forklinewidth depth \dp0
                        height \ht0 \hfil}
\setbox\noforkbox\hbox{\box1\box0\relax}
\def\unionstick{\mathop{\copy\noforkbox}\limits}
\def\nonfork#1#2_#3{#1\unionstick_{\textstyle #3}#2}
\def\nonforkin#1#2_#3^#4{#1\unionstick_{\textstyle #3}^{\textstyle
    #4}#2}
\setbox0\hbox{$\textstyle\bigcup$}
\setbox1\hbox to \wd0{\hfil{\sl /\/}\hfil}
\setbox2\hbox to \wd0{\hfil\vrule height \ht0 depth \dp0 width
                                \forklinewidth\hfil}
\newbox\doesforkbox
\setbox\doesforkbox\hbox{\box1\box0\relax}
\def\nunionstick{\mathop{\copy\doesforkbox}\limits}

\def\fork#1#2_#3{#1\nunionstick_{\textstyle #3}#2}
\def\forkin#1#2_#3^#4{#1\nunionstick_{\textstyle #3}^{\textstyle
    #4}#2}

\newcommand{\stickT}{%
\setbox255=\hbox{\raise1ex\hbox{$\hspace{0.2pt}\,\bullet\,$}}
\mathord{\rlap{\hbox to\wd255{\hss\hbox{$|$}\hss}}
\box255}
}
\newcommand{\stickS}{%
\setbox255=\hbox{\raise0.6ex\hbox{$\scriptstyle\bullet$}}
\mathord{\rlap{\hbox to\wd255{\hss\hbox{$\scriptstyle|$}\hss}}
\box255}
}

\author[S. Shelah]{Saharon Shelah}
\address{Einstein Institute of Mathematics,
The Hebrew University of Jerusalem,
9190401, Jerusalem, Israel; and\\
Department of Mathematics,
Rutgers University,
Piscataway, NJ 08854-8019, USA}
\urladdr{https://shelah.logic.at/}
\thanks{
Publication no. 288, summer `86. 
For versions up to 2019, the author thanks Alice Leonhardt for the beautiful typing. In the latest version, the author thanks an individual who wishes to remain anonymous for generously funding typing services, and thanks Matt Grimes for the careful and beautiful typing.
The author would like to thank the United
States -- Israel Binational Science Foundation for partially
supporting this research.
}

\makeatletter
\@namedef{subjclassname@2020}{\textup{2020} Mathematics Subject Classification}
\makeatother
\subjclass[2020]{FILL}
\date{September 8, 2023} 

\title[Strong Partition Relations Below the Power Set]{Strong Partition Relations Below the Power Set: Consistency --- Was Sierpinski Right? Vol. II}

\begin{document}
\makeatletter\def\shfiuwefootnote{\gdef\@thefnmark{}\@footnotetext}\makeatother\shfiuwefootnote{Version 2024-01-29. See \url{https://shelah.logic.at/papers/288/} for possible updates.}
\begin{abstract}
    We continue here \cite{Sh:276} (see the introduction there)
    but we do not rely on it.
    The motivation was a conjecture of Galvin stating that 
    $\big( 2^{\omega}\ge \omega_2 \big)$ + $\big( \omega_2\to [\omega_1]^{n}_{h(n)} \big)$ is consistent for a suitable $h : \omega \to \omega$. 
    In section 5 we disprove this and give similar negative results.
    In section 3 we prove the consistency of the conjecture replacing $\omega_2$ by $2^\omega$, which is quite large, starting
    with an Erd\H os cardinal. 
    In section 1 we present iteration lemmas which 
    are needed
    when we replace $\omega$ by a larger $\lambda$, and in section 4 we generalize a theorem of Halpern and Lauchli replacing $\omega$ by a larger $\lambda$. 
 
    This is a slightly corrected version of an old work.
\end{abstract}
\maketitle

\setcounter{section}{-1}
\section{Preliminaries}

Let $<^*_\chi$ be a well ordering of $\clH(\chi)$, where 
$$\clH(\chi) = \{x : \text{the transitive closure of $x$ has cardinality} < \chi\}$$ 
 agreeing with the usual well-ordering of the ordinals. 
$\bbP$ (and $\bbQ$, $\bbR$) will denote forcing notions; i.e. partial orders (really, quasiorders) with a minimal element $\varnothing = \varnothing_\bbP$. 

A forcing notion $\bbP$ is $\lambda$-closed if every increasing
sequence of members of $\bbP$ of length less than $\lambda$ has an
upper bound.

If $\bbP\in \clH(\chi)$, then for a sequence 
$\bar p = \LL p_i : i < \gamma\RR$ of members of $\bbP$, let
$$\name{\alpha} = \name{\alpha}_{\bar p} = \sup \! \big\{ \name{j} :
\{\beta_j : j < \name{j}\} \text{ has an upper bound in } \bbP\big\}$$ 
and define $\& \bar p$, the \emph{canonical upper bound} of 
$\bar{p}$, as follows:
\begin{enumerate}
    \item[(a)] It is the least upper bound of $\{p_i : i < \name{\alpha}\}$ in $\bbP$, if there exists such an element. 
\sn
    \item[(b)] If upper bounds of $\bar p$ exist but are not unique,  we choose the $<_\chi^*$-first upper bound. 
\sn
    \item[(c)] $p_0$, if (a) and (b) fail and $\gamma > 0$.
\sn
    \item[(d)] $\varnothing_\bbP$, if $\gamma = 0$.
\end{enumerate}
Let $p_0\;\&\; p_1$ be the canonical upper bound of 
$\LL p_\ell : \ell < 2\RR$.

\begin{notation}\label{x2}
1) Take $[a]^\kappa \defeq \{b\subseteq a : |b| = \kappa\}$ and
$[a]^{<\kappa} \defeq \bigcup\limits_{\theta<\kappa}[a]^\theta$.

\sn
2) For sets of ordinals $A$ and $B$, define $H^\OP_{A,B}$ as
the maximal order preserving bijection between initial segments of $A$
and $B$: i.e. it is the function with
domain $\{\alpha\in A : \otp(\alpha\cap A) < \otp(B)\}$ such that
$H^\OP_{A,B}(\alpha) = \beta$ \underline{iff} $\alpha\in A$, 
$\beta \in B$, and $\otp(\alpha \cap A) = \otp(\beta \cap B)$.
\end{notation}

\begin{definition}\label{x5}
$\lambda \to^+ (\alpha)_\mu^{<\aleph_0}$ holds if whenever $F$ is a 
function from $[\lambda]^{<\aleph_0}$ to $\mu$ and $C \subseteq \lambda$ 
is a club, then there is $A\subseteq C$ of order type $\alpha$ such 
that for any $w_1,w_2 \in [A]^{<\aleph_0}$, 
$|w_1| = |w_2| \Rightarrow F(w_1) = F(w_2)$.
\end{definition}

\begin{definition}\label{x11} 
$\lambda\to[\alpha]^n_{\kappa,\theta}$ 
if for every function $F$ from $[\lambda]^n$ to $\kappa$ there is 
$A \subseteq \lambda$ of order type $\alpha$ such that 
$\{F(w) : w \in [A]^n\}$ has 
cardinality
$\le\theta$.
\end{definition}

\begin{definition}\label{x14} 
A forcing notion $\bbP$ satisfies the Knaster condition (or 
`has property $K$') if for any $\{p_i : i < \omega_1\}\subset \bbP$ 
there is an uncountable $A\subset \omega_1$ such that the conditions
$p_i$ and $p_j$ are compatible whenever $i,j\in A$.
\end{definition}

\section{Introduction}

Concerning \ref{a2}--\ref{a8}, see Shelah \cite{Sh:80} and Shelah and Stanley \cite{Sh:154}, \cite{Sh:154a}.

\begin{definition}\label{a2} 
 A forcing notion $\bbQ$ satisfies $*^\eps_\mu$,
where $\eps$ is a limit ordinal $<\mu$, if Player \textbf{I}
has a winning strategy in the following game:
 
\sn 
 \textbf{Playing}: the play finishes after $\eps$ moves. In the $\alpha^\tthh$ move:
\begin{itemize}
    \item[Player \textbf{I} --] If $\alpha \ne 0$ he chooses 
    $\LL q_\zeta^\alpha : \zeta < \mu^+\RR$ such that 
    $q_\zeta^\alpha \in \bbQ$ and
    $$(\forall\beta<\alpha)(\forall \zeta<\mu^+) [p_\zeta^\beta \le q_\zeta^\alpha]$$ 
    and he chooses a regressive function $f_\alpha : \mu^+ \to \mu^+$ (i.e. $f_\alpha(i) < 1 + i$). 
    If $\alpha = 0$ let $q_\zeta^\alpha = \varnothing_\bbQ$ and $f_\alpha = \varnothing$.
\sn
    \item[Player \textbf{II} --] He chooses 
    $\LL p_\zeta^\alpha : \zeta < \mu^+\RR$ such that 
    $q_\zeta^\alpha \le p_\zeta^\alpha \in \bbQ$.
\end{itemize}

\sn
\textbf{The outcome}: Player \textbf{I} wins provided whenever
$\mu < \zeta < \xi < \mu^+$, $\cf(\zeta) = \cf(\xi) = \mu$, and
$\bigwedge\limits_{\beta<\eps} f_\beta(\zeta) = f_\beta(\xi)$, 
the set $\{p_\zeta^\alpha : \alpha < \eps\} \cup \{p_\xi^\alpha : \alpha < \eps\}$ has an upper bound in $\bbQ$.
\end{definition}

\begin{definition}\label{a5} 
We call $\big\LL \bbP_i,\,\bbQ_j : i \le i(*),\ j < i(*)\big\RR$ a 
$*_\mu^\eps$-\emph{iteration} provided that:
\begin{enumerate}
    \item[(a)] It is a $(<\mu)$-support iteration ($\mu$ is a regular cardinal).
\sn
    \item[(b)] If $i_1 < i_2 \le i(*)$ and $\cf(i_1) \ne \mu$ then $\bbP_{i_2} / \bbP_{i_1}$ satisfies $*_\mu^\eps$.
\end{enumerate}
\end{definition}

\begin{lemma}\label{a8}
If $\bfq = \LL \bbP_i,\,\bbQ_j : i \le i(*),\ j < i(*)\RR$ 
is a $(<\mu)$-support iteration and $(a)$ or $(b)$ or $(c)$ 
 below hold, then it is a $*_\mu^\eps$-iteration.
\begin{enumerate}[$(a)$]
    \item $i(*)$ is limit and $\bfq\rest j(*)$ is a 
    $*^\eps_\mu$-iteration for every $j(*) < i(*)$.
\sn
    \item $i(*) = j(*)+1$, $\bfq\rest j(*)$ is a 
    $*^\eps_\mu$-iteration, and $\bbQ_{j(*)}$ satisfies $*_\mu^\eps$ in $\bfV^{\bbP_{j(*)}}$.
\sn
    \item $i(*) = j(*)+1,$ $\cf(j(*)) = \mu^+$, $\bfq\rest j(*)$ is a $*^\eps_\mu$-iteration, and for every successor $i < j(*)$,
    $\bbP_{i(*)}/ \bbP_i$ satisfies $*^\eps_\mu$.
\end{enumerate}
\end{lemma}

\begin{PROOF}{\ref{a8}}
Left to the reader (after reading \cite{Sh:80} or \cite{Sh:154a}).
\end{PROOF}

\begin{theorem}\label{a11} 
Suppose $\mu = \mu^{<\mu} < \chi < \lambda$ and $\lambda$ is a 
strongly inaccessible $k_2^2$-Mahlo cardinal, where $k^2_2$ is a
suitable natural number (see \cite[3.6(2)]{Sh:289}), and assume 
$\bfV = \bfL$ for simplicity. 

\underline{Then} for some forcing notion $\bbP$:
\begin{enumerate}[$(A)$]
    \item $\bbP$ is $\mu$-complete, satisfies the $\mu^+$-c.c., has cardinality $\lambda$, and\\ $\bfV^\bbP \models ``2^\mu = \lambda"$.
\sn
    \item $\Vdash_\bbP \lambda\to [\mu^+]^2_3$ and even
    $\lambda\to[\mu^+]^2_{\kappa,2}$ for $\kappa < \mu$. 
\sn
    \item If $\mu = \aleph_0$ then $\Vdash ``\MA_\chi"$.
\sn
    \item If $\mu > \aleph_0$ \underline{then} $\Vdash_\bbP$ 
    ``for every $\mu$-complete forcing notion $\bbQ$ of cardinality $\le\chi$ satisfying $*_\mu^\eps$, and for any dense sets 
    $D_i \subseteq \bbQ$, for $i < i_0 < \lambda$, there is a 
    directed $G \subseteq \bbQ$ with 
    $\bigwedge\limits_i G \cap D_i \ne \varnothing$".
\end{enumerate}
\end{theorem}

As the 
    proof\footnote{In \cite{Sh:546}, full details are given for stronger theorems.} 
is very similar to \cite{Sh:276}
(particularly after reading section 3), we do not give details.
We shall define below only the systems needed to complete the proof. 
More general ones are implicit in \cite{Sh:289}.

\begin{convention}\label{a14} 
We fix a one to one function $\Cd = \Cd_{\lambda,\mu}$ from ${}^{\mu>}\lambda$ onto $\lambda$.
\end{convention}

\begin{remark}\label{a17}
Below we could have $\otp(B_x) = \mu^++1$ with little change.
\end{remark}

\begin{definition}\label{a20} 
Let $\mu < \chi < \kappa \le \lambda$, $\lambda = \lambda^{<\mu}$,
$\chi = \chi^{<\mu}$, $\mu = \mu^{<\mu}$.
\begin{enumerate}
    \item[1)] We call $x$ a $(\lambda,\kappa,\chi,\mu)$-\emph{pre-candidate} if $x = \LL a^x_u : u \in I_x \RR$, where for some set $B_x$ (unique, in fact):
    \begin{enumerate}
        \item[(i)] $I_x = [B_x]^{\leq 2}$ 
\sn
        \item[(ii)] $B_x$ is a subset of $\kappa$ of order type $\mu^+$.
\sn
        \item[(iii)] $a_u^x$ is a subset of $\lambda$ of cardinality
        $\le\chi$ closed under $\Cd$.
\sn
        \item[(iv)] $a_u^x \cap B_x = u$
\sn
        \item[(v)] $a_u^x \cap a_v^x \subseteq a^x_{u \cap v}$
\sn
        \item[(vi)] If $u,v \in I_x$ and $|u|=|v|$ then $a_u^x$ and $a_v^x$ have the same order type (and so $H^\OP_{a_u^x,a^x_v}$ maps $a_u^x$ onto $a_v^x$).
\sn
        \item[(vii)] If $u_\ell,v_\ell \in I_x$ and 
        $|u_\ell| = |v_\ell|$ for $\ell = 1,2$, 
        $|u_1\cup u_2| = |v_1\cup v_2|$, and 
        $H^\OP_{a^x_{u_1} \cup a^x_{u_2},a^x_{v_1} \cup a^x_{v_2}}$ 
        maps $u_\ell$ onto $v_\ell$ for $\ell=1,2$ then $H^\OP_{a^x_{u_1},a^x_{v_1}}$ and $H^\OP_{a^x_{u_2},a^x_{v_2}}$ are compatible. 
    \end{enumerate}

    \item[2)] We say $x$ is a $(\lambda,\kappa,\chi,\mu)$-candidate if it has the form $\LL M_u^x: u \in I_x\RR$, where
    \begin{enumerate}
        \item[$(\alpha)$]
        \begin{enumerate}
            \item[(i)] $\big\LL |M_u^x| : u \in I_x \big\RR$ is a $(\lambda,\kappa,\chi,\mu)$-precandidate (with $B_x$ defined as $\bigcup I_x$).
\sn
            \item[(ii)] $\tau_x$ is a vocabulary with $(\le \chi)$-many $(<\mu)$-ary place predicates and function symbols. 
\sn
            \item[(iii)] Each $M_u^x$ is a $\tau_x$-model.
\sn
            \item[(iv)] For $u,\,v\in I_x$ with $|u| = |v|$, $M_u^x\rest(|M_u^x|\cap|M_v^x|)$ is a model, and in fact an elementary submodel of $M_v^x$, $M_u^x$ and $M^x_{u\cap v}$.
        \end{enumerate}
\sn
        \item[$(\beta)$] For $u,\,v\in I_x$ with $|u| = |v|$, the function $H^\OP_{|M_u^x|,|M_v^x|}$ is an isomorphism from $M_u^x$ onto $M_v^x$.
    \end{enumerate} 

    \item[3)] We say the set $\gA$ is a $(\lambda,\kappa,\chi,\mu)$-\emph{system} if
    \begin{enumerate}
        \item[(A)] Each $x \in \gA$ is a $(\lambda,\kappa,\chi,\mu)$-candidate.
\sn
        \item[(B)] \textbf{Guessing}: if $\tau$ is as in (2)$(\alpha)$(ii) and $M^*$ is a $\tau$-model with universe $\lambda$, then for some $x\in\gA$, $s\in B_x \Rightarrow M_s^x \prec M^*$.
    \end{enumerate}
\end{enumerate} 
\end{definition}

\begin{definition}\label{a23} 1) We call the system $\gA$ \emph{disjoint} when:
\begin{itemize}
    \item[$(*)$] If $x\ne y$ are from $\gA$ and $\otp(|M_\varnothing^x|) \le \otp(|M_\varnothing^y|)$ \underline{then} for some 
    $B_1 \subseteq B_x$, $B_2 \subseteq B_y$ we have
    \begin{enumerate}[(a)]
        \item $|B_1| + |B_2| < \mu^+$
\sn
        \item The sets
        $$\bigcup\{|M_s^x|: s\in [B_x\setminus B_1]^{\le2}\} \text{ and } \bigcup\{|M_s^y|: s\in [B_y\setminus B_2]^{\le2}\}$$
        have intersection $\subseteq M_\varnothing^y$. 
    \end{enumerate}
\end{itemize}
\sn
2) We call the system $\gA$ \emph{almost disjoint} when:
\begin{itemize}
    \item[$(**)$] If $x,\,y\in \gA$ and 
    $\otp(|M_\varnothing^x|) \le \otp(|M_\varnothing^y|)$ \underline{then} for some $B_1\subseteq B_x$ and $B_2\subseteq B_y$ we have:
    \begin{enumerate}[(a)]
        \item $|B_1|+|B_2| < \mu^+$
\sn
        \item If $s\in [B_x\setminus B_1]^{\le2}$, 
        $t\in [B_y\setminus B_2]^{\le2}$ then $|M_s^x|\cap|M_t^x|\subseteq |M_\varnothing^y|$.
    \end{enumerate}
\end{itemize}
\end{definition}

\newpage

\section{Introducing the partition on trees}

\begin{definition}\label{b2} Let

\sn
1) $\Per({^{\mu>}2})$ be the set of $T$ such that
\begin{enumerate}
    \item $T\subseteq {}^{\mu>}2$, $\LL\ \RR\in T$.
\sn 
    \item $\big(\forall\eta\in T\big)\,\big(\forall\alpha < \lh(\eta) \big) \big[\eta\rest\alpha \in T \big]$
\sn
    \item If $\eta \in T \cap {^\alpha}2$ and $\alpha < \beta < \mu$ then for some $\nu \in T \cap {}^\beta2$ we have $\eta \lhd \nu$.
\sn 
    \item If $\eta \in T$ then for some $\nu$ we have $\eta \lhd \nu$, $\nu \caret \LL0\RR \in T$, and $\nu \caret \LL1\RR\in T$.
\sn 
    \item If $\eta \in {^\delta}2$, $\delta < \mu $ is a limit 
    ordinal, and $\{ \eta \rest \alpha : \alpha < \delta\} \subseteq T$ then $\eta \in T$.
\end{enumerate}

\sn
2) $\Per_\fe({^{\mu>}2})=$ 
$$\Big\{T\in\Per({}^{\mu>}2) : \alpha < \mu,\ \nu_1,\nu_2 \in {^\alpha2}\cap T \Rightarrow \Big[ \textstyle\bigwedge\limits_{\ell=0}^1 \nu_1 \caret \LL\ell\RR \in T \Leftrightarrow \textstyle\bigwedge\limits^1_{\ell=0} \nu_2 \caret \LL\ell\RR \in T \Big] \Big\}.$$

\noindent
3) $\Per_\uq({^{\mu>}2})=$ 
$$\Big\{T\in\Per({}^{\mu>}2 ): \alpha < \mu,\ \nu_1\ne\nu_2 \text{ from } {{}^\alpha}2\cap T\ \Rightarrow\ \textstyle\bigvee\limits^1_{\ell=0} 
 \textstyle\bigvee\limits^2_{m=1}\nu_m \caret \LL\ell\RR \notin T\Big\}$$

\noindent
4) For $T\in \Per({^{\mu>}2})$, let 
$\lim T = \big\{\eta \in {}^\mu2 : (\forall \alpha < \mu)[\eta\rest \alpha \in T]\big\}$.

\sn
5) For $T \in\Per_\fe({}^{\mu>}2)$ let $\clp_T : T \to {}^{\mu>}2$
be the unique one-to-one function from $\essp(T) \defeq \{\eta \in T : \eta \caret \LL0\RR, \eta \caret \LL1\RR \in T\}$ onto ${}^{\mu>}2$ which preserves $\lhd$ and lexicographic order.

\sn
6) Let $\SP(T) = \{\lh(\eta) : \eta \in \essp(T)\}$, and for 
$\eta,\nu \in T$ let
$$\essp(\eta,\nu) = \min\{i : \eta(i) \ne \nu(i) \vee i = \lh(\eta) \vee i = \lh(\nu)\}$$
(hence $\essp(\eta,\eta) = \lh(\eta)$).
\end{definition}

\begin{definition}\label{b5} 
For cardinals $\mu,\sigma$ and $n < \omega$ and 
$T \in \Per({}^{\mu>}2)$, let
\begin{enumerate}
    \item[1)] $\Col_\sigma^n(T) = \big\{d : d$ is a function from $\bigcup\limits_{\alpha<\mu}[{}^\alpha 2]^n\cap T$ to $\sigma\big\}$. We may write $d(\nu_0,\ldots,\nu_{n-1})$ for 
    $d \big( \{ \nu_0,\ldots,\nu_{n-1} \} \big)$.
\sn
    \item[2)] Let $<_\alpha^*$ denote a well ordering of ${}^\alpha2$ 
    (in this section it is arbitrary). We call $d \in \Col_\sigma^n(T)$ \emph{end-homogeneous} for $\big\LL {<_\alpha^*} : \alpha < \mu \big\RR$ provided that if $\alpha < \beta$ are from $\SP(T)$, 
    $\{\nu_0,\ldots,\nu_{n-1}\} \subseteq {}^\beta2\cap T$, $\LL\nu_\ell\rest\alpha : \ell < n\RR$ are pairwise distinct, and $\bigwedge\limits_{\ell,m} [\nu_\ell <_\beta^* \nu_m \Leftrightarrow \nu_\ell \rest \alpha <_\alpha^* \nu_m \rest \alpha]$ then
    $$d(\nu_0,\ldots,\nu_{n-1}) = d(\nu_0\rest\alpha,\ldots,\nu_{n-1} \rest \alpha).$$
    \item[3)] Let $\Eh\Col_\sigma^n(T) =$ 
    $$\big\{ d \in \Col_\sigma^n(T) : d \text{ is end-homogeneous for some } \LL {<_\alpha^*} : \alpha < \mu \RR\big\}$$
    ({see above}).
\sn
    \item[4)] For $\nu_0,\ldots,\nu_{n-1},\eta_0,\ldots,\eta_{n-1}$ 
    from ${}^{\mu>}2$, we say $\bar\nu = \LL\nu_0,\ldots,\nu_{n-1}\RR$ and $\bar\eta = \LL\eta_o, \ldots, \eta_{n-1}\RR$ are \emph{strongly similar} for $\big\LL {<_\alpha^*} : \alpha < \mu \big\RR$ if:
    \begin{enumerate}
        \item[(i)] $\lh(\nu_\ell)=\lh(\eta_\ell)$
\sn
        \item[(ii)] $\essp(\nu_\ell, \nu_m) = \essp(\eta_\ell, \eta_m)$ (equivalently, $ \lh(\nu_\ell \cap \nu_m) = \lh(\eta_\ell \cap \eta_m)$).
\sn
        \item[(iii)] If $\ell_1,\ell_2,\ell_3,\ell_4 < n$ and 
        $\alpha = \essp(\nu_{\ell_1},\,\nu_{\ell_2})$, $\alpha \leq \lh(\nu_{\ell_3}),\lh(\nu_{\ell_4})$, then
        $$\nu_{\ell_3} \rest \alpha <_{\alpha}^* \nu_{\ell_4} \rest \alpha \Leftrightarrow \eta_{\ell_3} \rest \alpha <_\alpha^* \eta_{\ell_4} \rest \alpha \text{ and } \nu_{\ell_3}(\alpha) = \eta_{\ell_3}(\alpha).$$
    \end{enumerate}
\sn
    \item[5)] For $\nu^a_0,\ldots,\nu^a_{n-1}$, 
    $\nu_0^b,\ldots,\nu_{n-1}^b$ from ${}^{\mu>}2$, we say 
    $\bar\nu^a = \LL\nu_0^a,\ldots,\nu^a_{n-1}\RR$ and 
    $\bar\nu^b = \LL\nu_0^b,\ldots,\nu^b_{n-1}\RR$ are \emph{similar} if the truth values of (i)--(iii) below do not depend on 
    $t \in \{a,b\}$ for any $\ell(1),\ell(2),\ell(3),\ell(4) < n$:
    \begin{enumerate}
        \item[(i)] $\lh(\nu_{\ell(1)}^t) < \lh(\nu_{\ell(2)}^t)$
\sn
        \item[(ii)] $\essp(\nu_{\ell(1)}^t,\,\nu_{\ell(2)}^t)< \essp(\nu_{\ell(3)}^t,\,\nu_{\ell(4)}^t)$
\sn
        \item[(iii)] for $\alpha = \essp(\nu_{\ell(1)}^t,\,\nu_{\ell(2)}^t)$ and 
        $ \lh(\nu_{\ell(3)}^t), \lh(\nu_{\ell(4)}^t) \geq \alpha$, the truth value of the following does not depend on $\ell$:
        $$\nu_{\ell(3)}^t \rest \alpha<_{\alpha}^* \nu_{\ell(4)}^t \rest \alpha \text{ and } \nu_{\ell(3)}^t(\alpha) = 0.$$
    \end{enumerate}
\sn
    \item[6)] We say $d\in\Col_\sigma^n(T)$ is almost homogeneous [homogeneous] on $T_1\subseteq T$ (for 
    $\big\LL {<_\alpha^*} : \alpha < \mu \big\RR$) \underline{if} for every $\alpha \in \SP(T_1)$, $\bar\nu$, 
    $\bar \eta \in [{^\alpha2}]^n \cap T_1$ which are strongly similar [similar] we have $d(\bar\nu) = d(\bar\eta)$.

    \item[7)] We say $\big\LL {<_\alpha^*} : \alpha < \mu\big\RR$ is 
    \emph{nice} to $T \in \Per({^{\mu>}2})$, provided that: if $\alpha<\beta$ are from $\SP(T)$, 
    $(\alpha,\beta) \cap \SP(T) = \varnothing$, 
    $\eta_1 \ne \eta_2 \in {}^\beta2 \cap T$, $\big[\eta_1 \rest \alpha <_\alpha^* \eta_2 \rest \alpha \hbox{ or } \eta_1 \rest \alpha = \eta_2 \rest \alpha, \> \eta_1(\alpha) < \eta_2(\alpha)\big]$ then $\eta_1 <_\beta^* \eta_2$. 
\end{enumerate}
\end{definition}




\begin{definition}\label{b8} 
1) $\Pr_\eht(\mu,n,\sigma)$ means ``for every $d \in \Col_\sigma^n({}^{\mu>}2)$, for some $T \in \Per({}^{\mu>}2)$,
$d$ is end homogeneous on $T$.''

\sn
2) $\Pr_\aht(\mu,n,\sigma)$ means ``for every
$d\in\Col_\sigma^n({}^{\mu>}2)$, for some $T\in\Per({}^{\mu>}2)$, 
$d$ is almost homogeneous on $T$.''

\sn
3) $\Pr_\hgt(\mu,n,\sigma)$ means for every
$d \in \Col_\sigma^n({}^{\mu>}2)$, for some $T\in\Per({}^{\mu>}2)$, 
$d$ is homogeneous on $T$.

\sn
4) For $x\in \{\eht,\aht,\hgt\}$, $\Pr_x^\fe(\mu,n,\sigma)$ 
is defined like $\Pr_x(\mu,n,\sigma)$ but we demand 
$T \in \Per_\fe({}^{\mu>}2)$.

\sn
5) If above we replace $\eht,\aht,\hgt$ by
$\ehtn,\ahtn,\htn$, respectively, this means\\
$\big\LL {<_\alpha^*} : \alpha < \mu \big\RR$ is fixed \emph{a priori}.

\sn
6) Replacing $n$ by ``$<\kappa$" and $\sigma$ by
$\olsi\sigma = \LL\sigma_\ell : \ell < \kappa\RR$ for 
$\kappa \le \aleph_0$ means that $\LL d_n : n < \kappa\RR$ are given, 
$d_n \in \Col_\sigma^n({}^{\mu>}2)$, and the conclusion holds for all 
$d_n$ with $n < \kappa$ simultaneously. Replacing ``$\sigma$" by 
``$<\sigma$" means that the assertion holds for every 
$\sigma_1 < \sigma$. 
\end{definition}

\begin{definition}\label{b11} 
1) $\Pr_\aht(\mu,n,\sigma(1),\sigma(2))$ means: for
every $d\in\Col_{\sigma(1)}^n({}^{\mu>}2)$, for some
$T\in\Per({}^{\mu>}2)$ and $\big\LL {<_\alpha^*} : \alpha < \mu \big\RR$, for every $\bar\eta \in \bigcup\big\{[{^{\alpha}}2]^n \cap
T : \alpha \in \SP(T)\big\}$, the set
$$
\Big\{d(\bar\nu) : \bar\nu \in \bigcup \big\{[{^\alpha}2]^n \cap
T_1  : \alpha \in \SP(T_1)\big\},\ \bar\eta
\text{ and $\bar\nu$ are strongly similar for }\big\LL {<_\alpha^*} : \alpha < \mu \big\RR \Big\}
$$
has cardinality $<\sigma(2)$.

\sn
2) $\Pr_\hgt(\mu,n,\sigma(1),\sigma(2))$ is defined similarly
with ``similar" instead of ``strongly similar".

\sn
3) $\Pr_x \! \big( \mu, <\kappa, \LL \sigma_\ell^1 : \ell < \kappa \RR, \> \LL \sigma_\ell^2 : \ell < \kappa \RR \big)$,
$\Pr_x^\fe(\mu,n,\sigma(1),\sigma(2))$,
$\Pr_x^\fe(\mu,<\aleph_0,\olsi\sigma^1,\olsi\sigma^2)$ are defined in the same way.
\end{definition}

There are many obvious implications.

\begin{fact}\label{b14} 
1) For every $T \in \Per({}^{\mu>}2)$ there is a
$T_1\subseteq T$ with $T_1\in \Per_\uq({}^{\mu>}2)$.

\sn
2) In defining $\Pr_x^\fe(\mu,n,\sigma)$ we can demand
$T \subseteq T_0$ for any $T_0\in \Per_\fe({}^{\mu>}2)$; similarly
for $\Pr_x^\fe(\mu,<\kappa,\sigma)$.

\sn
3) The obvious monotonicity holds.
\end{fact}

\begin{claim}\label{b17} 
1) Suppose $\mu$ is regular, 
$\sigma \ge \aleph_0$, and $\Pr^\fe_\eht(\mu,n,{<}\sigma)$. 
Then\\ $\Pr^\fe_\aht(\mu,n,{<}\sigma)$ holds.

\sn
1A) Similarly for $\Pr^\fe_\ehtn$ and $\Pr^\fe_\ahtn$.

\sn
2) If $\mu$ is weakly compact and
$\Pr^\fe_\aht(\mu,n,{<}\sigma)$ with $\sigma < \mu$, then 
$\Pr^\fe_\hgt(\mu,n,{<}\sigma)$ holds.
 
\sn
3) If $\mu$ is Ramsey and $\Pr^\fe_\aht(\mu,<\aleph_0,{<}\sigma)$ 
with $\sigma < \mu$, then
$\Pr^\fe_\hgt(\mu,<\aleph_0,{<}\sigma)$. 
 
\sn
4) If $\mu=\omega$, in the ``nice" version, the orders
$\big\LL {<_\alpha^*} :\alpha<\mu\big\RR $ disappear.

\sn
5) In parts (1)-(3), we can replace $\aht,\eht,\hgt$ by $\ahtn,\ehtn,\htn$ respectively.

\sn
6) In $\Pr_\eht^\fe(\mu,n,\sigma)$, we can strengthen the conclusion to:
\begin{enumerate}
    \item[$(*)$] If $\alpha < \beta$ are from $\SP(T)$, 
    $\LL \eta_\ell : \ell < n\RR \in {}^n(2^\alpha)$ is $<_\alpha$-increasing, and $\eta_\ell \lhd \nu_\ell^\iota \in 2^\beta$ 
    (for $\ell < n$ and $\iota \in \{1,2\}$) \underline{then}
    $$d(\{\nu_\ell^1 : \ell < n\}) = d(\{\nu_\ell^2 : \ell < n\}).$$
\end{enumerate}
\end{claim}

\begin{PROOF}{\ref{b17}} 
Easy; e.g. for (1A) we can use (6).

We induct on $n$; for $n+1$ and given 
$d_{n+1} : \bigcup \big\{ [{}^{\alpha}2]^{n+1} : \alpha < \mu \big\} \to \sigma$ 
and $\overbar <^{n+1} = \LL {<_\alpha^{n+1}} : \alpha < \mu \RR$, we apply $\Pr_\ehtn^\fe(\mu, n, {<\sigma})$. We get $T$. 

Let $f = \clp_T : T \to {}^{\mu>}2$ be as in \ref{b2}(5). 
Define $\overbar <^* = \LL {<_\alpha^*} : \alpha < \mu \RR$ and $d_n$ as follows:
\begin{enumerate}
    \item For $\alpha < \mu$ and $\eta_0,\eta_1 \in {}^{\alpha}2$, $\clp_T(\nu_\ell) = \eta_\ell$, $\lh(\nu_\ell) = \beta$ \underline{then} 
    $$\eta_0 <_\alpha^{n} \eta_1 \Leftrightarrow \nu_0 <_\alpha^{n+1} \nu_1$$

    \item for $\alpha < \mu$ and $\eta_0 <_\alpha^{n} \ldots <_\alpha^{n} \eta_{n-1}$, $\clp_T(\nu_\ell) = \eta_\ell$, $\lh(\nu_\ell) = \beta$ and for $k < n$, $\rho < 2$ we have $\nu_k \caret \LL \ell \RR \lhd \rho_{k,\ell} \in \essp(T_{n+1}) \cap {}^{\gamma}2$. {If} $\gamma$ is minimal then $d_n(\{ \eta_0, \ldots \eta_{n-1} \})$ codes the set of the following objects $\bft$:
    \begin{itemize}
        \item For some $\gamma > \alpha$ there are 
        $\rho_{k,\ell} \in \essp(T_{n+1}) \cap {}^{\gamma}2$ such that $\nu_k \caret \LL \ell \RR \unlhd \rho_{k,\ell}$ for $k < n$, $\ell < 2$ and $\bft$ codes all the information on the sequence $\LL \rho_{k,\ell} : k < n, \ell < 2 \RR$ (i.e. the order $<_\gamma^{n+1}$ and instances of $\bfd_{n+1}$). \qedhere
    \end{itemize}
\end{enumerate}
\end{PROOF}

The following theorem is a quite strong positive result for
$\mu = \omega$. Halpern-Lauchli proved \ref{b20}(1), Laver proved 
\ref{b20}(2) (and hence (3)), Pincus pointed out that Halpern-Lauchli's
proof can be modified to get \ref{b20}(2), and then
$\Pr^\fe_\eht(\omega,n,{<}\sigma)$ and (by it)
$\Pr^\fe_\hgt(\omega,n,{<}\sigma)$ are easy. 

\begin{theorem}\label{b20} 
$1)$ If $d\in\Col_\sigma^n({}^{\omega>}2)$ and $\sigma < \aleph_0$, 
\underline{then} there are $T_0,\ldots,T_{n-1} \in 
\Per_\fe({}^{\omega>}2)$ and
$k_0 < k_1 < \ldots < k_\ell < \ldots$ and $s < \sigma$ such that for every $\ell < \omega$, \underline{if} $\mu_0 \in T_0$, $\mu_1 \in T_1,\ldots,\nu_{n-1} \in
T_{n-1},\; \bigwedge\limits_{m<n}\lh(\nu_m)=k_\ell$, then
$d(\nu_0,\ldots,\nu_{n-1})=s$. 

\sn
$2)$ We can demand in 1) that 
$$\SP(T_\ell) = \{k_0,k_1,\ldots\}.$$

\sn
$3)$ $\Pr_\htn^\fe(\omega,n,\sigma)$ for $\sigma < \aleph_0$.

\sn
$4)$ $\Pr_\htn^\fe\big(\omega,<\aleph_0,\LL\sigma_n^1 : n < \omega \RR,
\LL \sigma_n^2 : n < \omega \RR\big) $ if $\sigma_n^1 < \aleph_0$ and
$\LL\sigma_n^2 : n < \omega \RR$ diverge to infinity.
\end{theorem}

\begin{definition}\label{b23} 
Let $d$ be a function with domain $\supseteq [A]^n$, 
$A$ be a set of ordinals, $F$ be a one-to-one 
function from $A$ to ${}^{\alpha(*)}2$, $<_\alpha^*$ be a well 
ordering of ${}^\alpha2$ for $\alpha \le \alpha(*)$ such that 
$F(\alpha)<_\alpha^* F(\beta) \Leftrightarrow \alpha<\beta$, and
$\sigma$ be a cardinal.

\sn
1) We say $d$ is $(F,\sigma)$-\emph{canonical} on $A$ if for any
$\alpha_1 < \ldots < \alpha_n \in A$,
$$\left\vert \big\{ d(\beta_1,\ldots,\beta_n) : \big\LL F(\beta_1),\ldots,F(\beta_n) \big\RR \text{ similar to } \big\LL F(\alpha_1),\ldots,F(\alpha_n) \big\RR \big\}\right\vert\le\sigma$$

\sn
2) We define ``almost $(F,\sigma)$-canonical" similarly using
strongly similar instead of ``similar".
\end{definition}

\newpage

\section{Consistency of a strong partition below the continuum} 

This section is dedicated to the proof of

\begin{theorem}\label{c2} 
Suppose $\lambda$ is the first Erd\H os cardinal (i.e. the first such 
that\\ $\lambda \to (\omega_1)^{<\omega}_2$).
\underline{Then}, if $A$ is a Cohen subset of $\lambda$, in $\bfV[A]$ for
some $\aleph_1$--c.c. forcing notion $\bbP$ of cardinality $\lambda$,
$\Vdash_\bbP$ ``$\MA_{\aleph_1}(Knaster) + 2^{\aleph_0} = \lambda$" and:
\begin{enumerate}[1)]
    \item $\Vdash_\bbP ``\lambda\to [\aleph_1]^n_{h(n)}$" for suitable $h : \omega \to \omega$ (explicitly defined below).
\sn
    \item In $\bfV^\bbP$, for any colorings $d_n$ of $\lambda$ where $d_n$ is $n$-place, and for any divergent $\LL \sigma_n : n < \omega \RR$ (see below), there is a $W \subseteq \lambda$, $|W| = \aleph_1$ and a function $F : W \to {}^\omega2$ such that $d_n$ is $(F,\sigma_n)$-canonical on $W$ for each $n$. (See Definition \ref{b23} above.)
\end{enumerate}
\end{theorem}

\begin{remark}\label{c8} 
1) $h(n)$ is $n!$ times the number of $u\in[{^\omega2}]^n$ satisfying ``if $\eta_1,\eta_2,\eta_3,\eta_4\in u$ are distinct and $\eta_1 \cap \eta_2 \neq \eta_3 \cap \eta_4$ \underline{then} $\essp(\eta_1,\eta_2),$ $\essp(\eta_3,\eta_4)$ are distinct'' up to strong similarity for any nice 
$\big\LL{<^*_\alpha} : \alpha < \omega \big\RR.$

\sn
2) A sequence $\LL\sigma_n : n < \omega \RR$ is \emph{divergent} if
$(\forall m) (\exists k)(\forall n \ge k)[\sigma_n \ge m]$.
\end{remark}

\begin{notation}\label{c11} 
For a sequence $a = \LL a_i,\,e_i^* : i < \alpha\RR$ with 
$a_i \subseteq i$ and $e_i \in \{1,2\}$,
we call $b \subseteq \alpha$ \emph{closed} (or `$a$-{closed}') if
\begin{enumerate}
    \item[(i)]  $i\in b \Rightarrow a_i \subseteq b$
\sn
    \item[(ii)] If $i < \alpha,$ $e_i^*=1$, and $\sup(b \cap i) = i$ then $i \in b$.
\end{enumerate}
\end{notation}

\begin{definition}\label{c14} 
Let $\gK$ be the family of 
$\bfq = \LL \bbP_i,\Name\bbQ_j,a_j,e_j^* : j < \alpha,\ i \le \alpha\RR$ such that:
\begin{enumerate}
    \item[(a)] $a_i\subseteq i$, $|a_i|\le\aleph_1$, and 
    $e_i^* \in \{0,1\}$.
\sn
    \item[(b)] $a_i$ is closed for $\LL a_j,e_j^* : j < i\RR$
    and $[e_i^* = 1 \Rightarrow \cf(i) = \aleph_1]$.
\sn
    \item[(c)] $\bbP_i$ is a forcing notion, $\Name\bbQ_j$ is a $\bbP_j$-name of a forcing notion of cardinality $\aleph_1$ with minimal element $\varnothing$ or $\varnothing_j$, and for simplicity the underlying set of $\Name\bbQ_j$ is $\subseteq[\omega_1]^{<\aleph_0}$ (we do not lose {anything} by this).
\sn
    \item[(d)] $\bbP_\beta = \big\{p : p$ is a function whose domain is a finite subset of $\beta$ and for $i\in \dom(p)$,
    $\Vdash_{\bbP_i} ``f(i) \in \Name\bbQ_i" \big\}$ 
    with the order $p\le q$ if and only if for $i\in \dom(p)$, 
    $q\rest i \> \Vdash_{\bbP_i} ``p(i)\le q(i)$".
\sn
    \item[(e)] For $j<i$, $\Name\bbQ_j$ is a $\bbP_j$-name involving 
    only antichains contained in\\ $\{p\in \bbP_j : \dom(p) \subseteq a_j\}$.
\end{enumerate}
\textbf{Notation:} For $p\in \bbP_i$, $j<i$, $j \notin \dom(p)$ we let $p(j) = \varnothing$. Note that for $p\in \bbP_i$ {and} $j\le i$, {we have} $p\rest j \in \bbP_j$.  
\end{definition}

\begin{definition}\label{c17} 
For $\bfq \in \gK$ as above (so $\alpha = \lh(\bfq)$):

\sn
1) for any $b \subseteq \beta \le \alpha$ closed for 
$\LL a_i,e_i^* : i < \beta \RR$, we define $\bbP^\cn_b$ [by simultaneous induction on $\beta$]: 
$$
    \bbP^\cn_b = \left\{ p\in \bbP_\beta : \dom(p) \subseteq b, 
    \hbox{ and for } i \in \dom(p),\ p(i) \hbox{ is a canonical name}\right\}.
$$
I.e. for any $x$, $\{p \in \bbP^\cn_{a_i} : p \Vdash_{\bbP_i} ``p(i) = x$''
or $p\Vdash_{\bbP_i}``p(i)\ne x"\}$ is a predense subset of $\bbP_i$.
 
\sn
2) For $\bfq $ as above, $\alpha = \lh(\bfq)$, take $\bfq \rest \beta = \LL \bbP_i, \Name \bbQ_j,a_j : i \le \beta,\ j < \beta \RR$
for $\beta \le \alpha$ and the order is the order in $\bbP_\alpha$
(if $\beta\ge \alpha$, $\bfq\rest\beta = \bfq$).

\sn
3) ``\emph{$b$ closed for $\bfq$} means ``$b$ closed for 
$\LL a_i,e^*_i : i < \lh(\bfq)\RR$''.
\end{definition}

\begin{fact}\label{c20} 
1) if $\bfq\in \gK$ then $\bfq\rest\beta\in\gK$.

\sn
2) Suppose $b \subseteq c \subseteq \beta \le \lh(\bar\theta)$,
$b$ and $c$ are closed for $\bfq \in \gK$. 
\begin{enumerate}
    \item[(i)] If $p \in \bbP^\cn_c$ then $p \rest b \in \bbP^\cn_b$.
\sn 
    \item[(ii)] If $p,q \in \bbP_c^\cn$ and $p \le q$ then 
    $p \rest b \le q \rest c$.
\sn 
    \item[(iii)] $\bbP_c^\cn \lessdot \bbP_\beta$.
\end{enumerate}

\sn
3) $\lh(\bfq)$ is closed for $\bfq$.

\sn
4) If $\bfq \in \gK$, $\alpha = \lh(\bfq)$ then
$\bbP^\cn_\alpha$ is a dense subset of $\bbP_\alpha$.

\sn
5) If $b$ is closed for $\bfq$, $p,q\in \bbP^\cn_{\lh(\bfq)}$,
$p\le q$ in $\bbP_{\lh(\bfq)}$ and $i \in \dom(p)$ 
then\\ $q \rest a_i \Vdash_{\bbP_i} ``p(i)\le q(i)$''
hence $\Vdash_{\bbP^\cn_{a_i}} ``p(i)\le_{\bbQ_i} q(i)$''.
\end{fact}

\begin{definition}\label{c23} 
Suppose $W = (W,\le)$ is a finite partial order and 
$\bfq \in \gK$.

\sn
1) $\IN_W(\bfq)$ is the set of $\bar b$-s
satisfying $(\alpha)$--$(\gamma)$ below:
\begin{enumerate}
    \item[$(\alpha)$] $\bar b = \LL b_w : w \in W\RR$ is an 
    indexed set of $\bfq$-closed subsets of $\lh(\bfq)$.
\sn
    \item[$(\beta)$] $W \models w_1 \le w_2 \Rightarrow b_{w_1} \subseteq b_{w_2}$.
\sn
    \item[$(\gamma)$] {If} $\zeta \in b_{w_1} \cap b_{w_2}$, 
    $w_1 \le w$, and $w_2 \le w$ then 
    $$(\exists u\in W)[\zeta \in b_u \wedge u \le w_1 \wedge u \le w_2].$$
\end{enumerate}
We assume $\bar b$ codes $(W,\le)$.

\sn
2) For $\bar b\in \IN_W(\bfq)$, let
$$\bfq[\bar b] = \big\{\LL p_w : w \in W\RR : p_w \in
P^\cn_{b_w},\ [W \models w_1 \le w_2\Rightarrow
p_{w_2}\rest b_{w_1} = p_{w_1}]\big\}$$ with ordering
$\bfq[\bar b] \models \bar p^1 \le \bar p^2$ iff
$\bigwedge\limits_{w\in W} p^1_w \le p^2_w$.

\sn
3) Let $\gK^1$ be the family of $\bfq \in \gK$ such 
that for every $\beta \le \lh(\bfq)$ and $(\bfq\rest\beta )$-closed set $b$, $\bbP_\beta$ and $\bbP_\beta / \bbP^\cn_b$ satisfy the Knaster condition.
\end{definition}

\begin{fact}\label{c26} 
Suppose $\bfq\in \gK^1$, $(W,\le)$ is a finite partial order, 
$\bar b \in \IN_W(\bfq)$ and $\bar p \in\bfq[\bar b]$.

\sn
1) If $w\in W$, $p_w \le q \in \bbP_{b_w}^\cn$ \underline{then} there
is $\bar r \in\bfq[\bar b]$, $q \le r_w$, $\bar p \le \bar r$. In fact, 
\[
    r_u(\gamma ) = 
    \begin{cases}
        p_u(\gamma )& \text{if } \gamma \in \dom\ p_u \setminus \dom\ q,\\ 
        p_u(\gamma )\,\, \&\,\, q(\gamma )& \text{if } \gamma\in b_u\cap \dom\ q \text{ and for some } v \in W,\\ 
&\quad u \leq v \le w \text{ and } \gamma\in b_v,\cr
p_u(\gamma)& \text{if } \gamma\in b_u\cap\dom\ q 
        \text{ but the previous case fails.} 
    \end{cases}
\]

\noindent
2) Suppose $(W_1,\le )$ is a submodel of $(W_2,\le )$, both
finite partial orders, $\bar b^l\in \IN_{W_l}(\bfq)$, 
$\bar b^1_w =\bar b^2_w$ for $w\in W_1$. 
\begin{enumerate}
    \item[$(\alpha)$] If $\bar q\in\bfq[\bar b^2]$ then 
    $\LL q_w : w \in W_1\RR \in\bfq[\bar b^1]$.
\sn
    \item[$(\beta)$] If $\bar p \in\bfq[\bar b^1]$ \underline{then} there is $\bar q\in\bfq[\bar b^2]$ with
    $\bar q\rest W_1 = \bar p$; in fact, $q_w(\gamma)$ is $p_u(\gamma)$ if $u \in W_1$, $\gamma \in b_u$,and $u \le w$, provided that
    \begin{enumerate}
        \item[$(**)$] If $w_1,w_2\in W_1$, $w\in W_2$,
   $w_1\le w$, $w_2\le w$ and $\zeta \in b_{w_1} \cap b_{w_2}$
        then for some $v\in W_1$, $\zeta\in b_v$, $v\le w_1$, $v\le w_2$.
    \end{enumerate}
    (This guarantees that if there are several $u$-s as above we shall get the same value.)
\end{enumerate}

\sn
3) If $\bfq\in\gK^1$ then $\bfq[\bar b]$ satisfies the 
Knaster condition. If $\varnothing$ is the minimal element of $W$ 
(i.e. $u\in W \Rightarrow W\models \varnothing\le u$)
then $\bfq[\bar b]/ \bbP^\cn_{b_\varnothing}$ also satisfies the Knaster
condition and so {is} $\lessdot\ \bfq[\bar b]$, when we
identify $p\in \bbP^\cn_b$ with $\LL p : w \in W \RR$.
\end{fact}

\begin{PROOF}{\ref{c26}} 
1) It is easy to check that each $r_u(\gamma)$ is in
$\bbP^\cn_{b_u}$. So, in order to prove $\bar r\in\bfq[\bar b]$, 
we assume $W\models u_1\le u_2$ and have to prove that 
$r_{u_2} \rest b_{u_1} = r_{u_1}$. Let $\zeta\in b_{u_1}$.

\mn
\textbf{First case:} $\zeta \notin \dom(p_{u_1}) \cup \dom (q)$. 

So $\zeta \notin \dom(r_{u_1})$ (by the definition of $r_{u_1}$) and
$\zeta \notin \dom (p_{u_2})$ (as $\bar p \in \bfq[\bar b]$) hence
$\zeta \notin \dom (p_{u_2}) \cup \dom(q)$ hence $\zeta \notin
\dom(r_{u_2})$ by the choice of $r_{u_2}$, so we have finished.

\mn
\textbf{Second case:} $\zeta\in\dom (p_{u_1})\setminus \dom(q)$.

As $\bar p\in\bfq[\bar b]$ we have $p_{u_1}(\zeta)=p_{u_2}(\zeta)$,
and by their definition, $r_{u_1}(\zeta)=p_{u_1}(\zeta)$,
$r_{u_2}(\zeta)=p_{u_2}(\zeta)$.

\mn
\textbf{Third case:} $\zeta\in\dom(q)$ and $(\exists v\in W)$
$[\zeta\in b_v \land v\le u_1\land v\le w]$. 
 
By the definition of $r_{u_1}(\zeta)$, we have 
$r_{u_1}(\zeta)= p_{u_1}(\zeta)\ \&\ q(\zeta)$; also,
the same $v$ witnesses $r_{u_2}(\zeta)=p_{u_2}(\zeta)\ \&\ q(\zeta)$ 
$$(\text{as } \zeta \in b_v \land v \le u_1 \land v \le w \Rightarrow 
 {\zeta \in b_v} \land v \le u_2 \land 
{v \le w}),$$ 
and of course $p_{u_1}(\zeta) = p_{u_2}(\zeta)$ (as
$\bar p\in\bfq[\bar b]$).

\mn
\textbf{Fourth case:} $\zeta\in\dom(q)$ and $\neg(\exists v\in
W)[\zeta\in b_v\land v\le u_1 \land v\le w]$. 

By the definition of $r_{u_1}(\zeta)$ we have
$r_{u_1}(\zeta) = p_{u_1}(\zeta)$. It is enough to prove that
$r_{u_2}(\zeta) = p_{u_2}(\zeta)$ as we know that
$p_{u_1}(\zeta) = p_{u_2}(\zeta)$ (because 
$\bar p \in\bfq[\bar b]$, $u_1 \le u_2)$. If not, then 
for some $v_0 \in W$, $\zeta \in b_{v_0}\land v_0\le u_2 \land v_0\le w$. But $\bar b\in
\IN_W(\bfq)$, hence (see condition
$(\gamma)$ of Definition \ref{c23}(1), applied with $\zeta,\,w_1,\,w_2,\,w$ there standing for
$\zeta,\,v_0,\, u_1,\,u_2$ here) we know that for some $v\in W$,
$\zeta\in v\land v\le v_0\land v\le u_1$. As $(W,\le)$ is a
partial order, $v\le v_0$ and $v_0\le w$, we can conclude $v\le
w$. So $v$ contradicts our being in the fourth case. So we have
finished the fourth case. 

Hence we have finished proving $\bar r\in\bfq[\bar b]$. We
also have to prove $q\le r_w$, but for $\zeta\in\dom(q)$ we have
$\zeta\in b_w$ (as $q\in \bbP^\cn_{w}$ is on assumption) and
$r_w(\zeta)=q(\zeta)$ because $r_w(\zeta)$ is defined
by the second case of the definition as 
    $$(\exists v\in W)[\zeta\in b_w\land v\le w \land v\ge w]$$ 
    i.e. $v=w$. 

Lastly, we have to prove that $\bar p\le \bar r$ (in $\bfq[\bar
b]$). So let $u\in W$, $\zeta \in \dom (p_u)$ and we have to prove
$r_u\rest \zeta\Vdash_{\bbP_\zeta}$``$p_u(\zeta)\le_{\bbP_\zeta} r_u(\zeta)$". As
$r_u(\zeta)$ is $p_u(\zeta)$ or $p_u(\zeta)\; \&\; q(\zeta)$ this is obvious. 

\mn
2) Immediate.

\mn
3) 
We prove this by induction on $|W|$.

For $|W| = 0$ this is totally trivial.

For $|W| = 1,2$ this is assumed.

For $|W| > 2$ fix $\bar p^i \in\bfq[\bar b]$ for $i < \omega_1$.
Choose a maximal element $v \in W$ and let 
$c = \bigcup \{b_w : W \models w < v\}$. 
Clearly $c$ is closed for $\bfq$. 

We know that $\bbP_c^\cn$, $\bbP^\cn_{b_v}/ \bbP_c^\cn$ are Knaster
by the induction hypothesis.
We also know that $p^i_v\rest c\in \bbP_c^\cn$ for $i < \omega_1$, hence
for some $r\in \bbP_c^\cn$,
$$r\Vdash ``\Name{A} = \left\{ i < \omega_1 : p^i_v \rest c \in \Name{G}_{\bbP_c^\cn} \right\} \hbox{ is
uncountable"}$$ 
hence
\begin{align*}
\Vdash &\text{``there is an uncountable } A^1 \subseteq
\Name{A}\text{ such that }\\ 
 &\ \left[ i,j \in A^1 \Rightarrow p^i_v,\,p_v^j \text{ are compatible in } \bbP^\cn_{b_v}/\Name{G}_{\bbP_c^\cn} \right]".
\end{align*}
Fix a $\bbP_c^\cn$-name $\Name{A}^1$ for such an $A^1$.

Let $A^2 = \left\{ i < \omega_1 : (\exists q \in P_c^\cn) [q\Vdash i \in \Name{A}^1]\right\}$. Necessarily 
$|A_2| = \aleph_1$, and for $i\in A^2$ there is $q^i\in \bbP_c^\cn$,
$q^i\Vdash i \in A^1$, and
without loss of generality~$p^i_v\rest c \le q^i$. Note that $p^i_v\ \&\ q^i\in \bbP^\cn_c$. 

For $i\in A^2$, let $\bar r^i$ be defined using \ref{c26}(1) (with
$\bar p^i$, $p^i_v\, \&\, q^i$). 
Let $W_1 = W \setminus \{v\}$, $\bar b' = \LL b_w : w \in W_1\RR$.

By the induction hypothesis applied to $W_1$, $\bar b'$,
$\bar r^i\rest W_1$, for $i\in A^2$ there
is an uncountable $A^3 \subseteq A^2$ and for $i < j$ in $A^3$, there is
$\bar r^{i,j}\in\bfq[\bar b']$ with $\bar r^i\rest W_1 \le \bar
r^{i,j}$ and $\bar r^j\rest W_1 \le \bar r^{i,j}$. Now 
define $r_c^{i,j} \in \bbP_c^\cn$ as follows: its domain is 
$\bigcup\left\{\dom(r_w^{i,j}) : W \models w < v\right\}$ 
and
$r_c^{i,j}\rest \dom(r_w^{i,j}) = r_w^{i,j}$ whenever $W \models w < v$. 

Why is this a definition? As $W \models w_1 \le v \wedge w_2 \le v$, 
$\zeta\in b_{w_1} \wedge \zeta \in b_{w_2}$ implies that for some 
$u \in W$, $u \le w_1\wedge u \le w_2$ and $\zeta \in u$. It is easy to
check that $r_c^{i,j}\in \bbP_c^\cn$. Now
$r_c^{i,j}\Vdash_{\bbP_c^\cn}$ ``$p^i_{b_v}, p^j_{b_v}$ are
compatible in $\bbP^\cn_{b_v}/ \bbP_c^\cn$". 

So there is $r\in \bbP^\cn_{b_v}$ such that $r^{i,j}_c \le r$,
$p^i_{b_v} \le r$, $p^j_{b_v} \le r$. 
As in part (1) of \ref{c26}, we can combine $r$ and $\bar r^{i,j}$ to a common
upper bound of $\bar p^i$, $\bar p^j$ in $\bfq[\bar b]$. 
\end{PROOF}

\begin{claim}\label{c29} 
If $e = 0,1$ and $\delta$ is a limit ordinal, and $\bbP_i,\Name\bbQ_i$,
$\alpha_i,e^*_i$ (for $i < \delta$) are such that for each $\alpha < \delta$, $\bfq^\alpha = \LL \bbP_i,\Name\bbQ_j,\alpha_j,e^*_j : i \le \alpha,\,\,j < \alpha \RR$ belongs to $\gK^\ell$, then for a unique $\bbP_\delta$, 
$\bfq = \LL \bbP_i, \Name\bbQ_j, \alpha_j, e^*_j : i \le \delta,\ j < \delta \RR$ belongs to $\gK^\ell$.
\end{claim}

\begin{PROOF}{\ref{c29}} 
We define $\bbP_\delta$ by Definition \ref{c14}(d). 
The least easy problem is to verify the Knaster conditions 
(for $\bfq\in \gK^1$). The proof is like the preservation 
of the c.c.c. under iteration for limit stages. 
\end{PROOF}

\begin{convention}\label{c32}
In \ref{c29}, we shall not make a strict distinction between\\ 
$\LL \bbP_i,\Name\bbQ_j,\alpha_j,e^*_j : i \le \delta,\ j < \delta\RR$ and $\LL \bbP_i,\Name\bbQ_i,\alpha_i,e^*_i : i < \delta\RR$.
\end{convention}

\begin{claim}\label{c35} 
If $\bfq \in \gK^\ell$, $\alpha = \lh(\bfq)$,
$a \subset \alpha$ is closed for $\bfq$, $|a| \le \aleph_1$, and
$\Name\bbQ_1$ is a $\bbP^\cn_a$-name of a forcing 
notion satisfying (in $\bfV^{\bbP_\alpha}$) the Knaster condition
whose underlying set is a subset of $[\omega_1]^{<\aleph_0}$,
\underline{then} there is a unique $\bfq^1\in \gK^\ell$ with 
$\lh(\bfq_1) = \alpha + 1$, $\bbQ^1_\alpha = \Name\bbQ$, and
$\bfq \rest \alpha =\bfq$.
\end{claim}

\begin{PROOF}{\ref{c35}}
Left to the reader.
\end{PROOF}

\bn
{We are now ready to prove \ref{c2}.}

\begin{PROOF}{\ref{c2}}
\textbf{Stage A:} We force by $\gK^1_{<\lambda} = \left\{\bfq \in \gK^1 : \lh(\bfq) < \lambda,\ \bfq \in \clH(\lambda )\right\}$
ordered by being an initial segment (which is 
equivalent to forcing a Cohen subset of $\lambda$). The generic
object is essentially $\bfq^*\in \gK^1_\lambda$, $\lh(\bfq^*) = \lambda$, and then we force by $\bbP_\lambda = \lim\bfq^*$. Clearly $\gK^\ell_{<\lambda}$ is a $\lambda$-complete forcing notion of cardinality $\lambda$, 
and $\bbP_\lambda$ satisfies the c.c.c. 
Clearly it suffices to prove part (2) of \ref{c2}.

Suppose $\name d_n$ is a name of a function from 
$[\lambda]^n$ to $\name k_n$ for $n < \omega$, 
$\name\sigma_n < \omega$,\\ $\LL \sigma_n : n < \omega \RR$ 
    diverges\footnote{I.e. $(\forall m)  (\exists k)  (\forall n\ge k)  [\sigma_n\ge m]$.}
and for some $\bfq^0 \in \gK^1_{<\lambda}$, we have
\begin{align*}
    \bfq^0 \Vdash _{\gK^1_{<\lambda}}  (\exists p\in \Name \bbP_\lambda) &
    \big[p\Vdash_{\bbP_\lambda} ``\LL \name d_n : n < \omega \RR \text{ is a}\\
    & \text{ counterexample to \ref{c2}(2)"}\big].
\end{align*}

In $\bfV$ we can define $\LL\bfq^\zeta : \zeta < \lambda \RR$ with
$\bfq^\zeta \in \gK^1_{<\lambda}$ such that
$$\zeta < \xi \Rightarrow \bfq^\zeta = \bfq^\xi \rest \lh(\bfq^\zeta).$$ 
In
$\bfq^{\zeta +1}$, $e^*_{\lh(\bfq_\zeta)} = 1$, $\bfq^{\zeta +1}$
forces (in $\gK^1_{<\lambda})$ a value to $p$ and the 
$\Name \bbP_\lambda$-names $\name d_n\rest\zeta$, $\name\sigma_n$,
$\name k_n$ for $n < \omega$; i.e. the values here are still 
$\bbP_\lambda$-names. Let $\bfq^*$ be the limit of the $\bfq^\xi$-s. So 
$\bfq^* \in \gK^1$, $\lh(\bfq^*) = \lambda$,
$\bfq^* = \LL \bbP^*_i,\Name\bbQ^*_j,\alpha^*_j,e^*_j : i \le \lambda,\ 
j < \lambda\RR$, and the $\bbP^*_\lambda$-names $\name d_n$,
$\name\sigma_n$, $\name k_n$ are defined such that in $\bfV^{\bbP^*_\lambda}$,
$\name d_n$, $\name\sigma_n$, $\name k_n$ contradict clause (2)
(as any $\bbP^*_\lambda$-name of a bounded subset of $\lambda$
is a $\bbP^*_{\lh(\bfq^\xi)}$-name for some $\xi < \lambda$). 
 
\mn
\textbf{Stage B:} 
Let $\chi = \kappa^+$ and $<^*_{\chi}$ be a well-ordering of $\clH(\chi)$.
Now we can apply $\lambda \to
(\omega_1)^{<\omega}_2$ to get $\delta,B,N_s$ and
$\bfh_{s,t}$ (for $s,t \in [B]^{<\aleph_0}$ with $|s| = |t|)$ such
that: 
\begin{enumerate}
\item[(a)] $B \subseteq \lambda$ with $\otp(B) = \omega_1$ and 
$\sup B = \delta$.
\sn
\item[(b)] $N_s \prec (\clH(\chi),\in,<^*_\chi)$, 
$\bfq^* \in N_s$,
$\LL \name d_{  
{n}},\name\sigma_n,\name k_n : n < \omega \RR \in N_s$.
\sn 
\item[(c)] $N_s \cap N_t = N_{s\cap t}$
\sn 
\item[(d)] $N_s \cap B = s$
\sn 
\item[(e)] If $s = t\cap \alpha$, $t\in[B]^{<\aleph_0}$ 
then $N_s \cap \lambda$ is an initial segment of $N_t$. 
\sn 
\item[(f)] $\bfh_{s,t}$ is an isomorphism
from $N_t$ onto $N_s$ (when defined). 
\sn 
\item[(g)] $\bfh_{t,s} = \bfh_{s,t}^{-1}$
\sn 
\item[(h)] $p_0\in N_s$, $p_0 \Vdash_{\bbP_\lambda}$ ``$\LL \name
d_n,\name\sigma_n,\name k_n : n < 
{\omega}\RR$ is a counterexample to the conclusion of \ref{c2}''.
\sn 
\item[(i)] $\omega_1 \subseteq N_s$, $|N_s| = \aleph_1$ and if
$\gamma \in N_s$, $\cf(\gamma) > \aleph_1$ then 
$\cf(\sup(\gamma \cap N_s)) = \omega_1$. 
\end{enumerate}

Let $\bfq = \bfq^*\rest\delta$, $\bbP = \bbP^*_{\delta}$
and $\bbP_a = \bbP^\cn_a$ (for $\bfq$), where $a$ is closed for $\bfq$.

\noindent 
Note: $\bbP^*_\lambda\cap N_s = \bbP^*_\delta \cap N_s = \bbP_{\sup \lambda\cap N_s} \cap N_s = \bbP_s\cap N_s$. Note also 
$$\gamma \in \lambda \cap N_s \Rightarrow a^*_\gamma \subseteq \lambda \cap N_s.$$

\mn
\textbf{Stage C:} It suffices to show that we can define
$\Name\bbQ_\delta$ in $\bfV^{\bbP_\delta}$ which forces a subset $W$ 
of $B$ of cardinality $\aleph_1$ and an $\Name F : W \to {}^\omega2$ which 
exemplify the desired conclusion in $(2)$, and prove that 
$\Name\bbQ_\delta$ satisfies the $\aleph_1$-c.c.c. in $\bfV^{\bbP_\delta}$ (and has cardinality $\aleph_1$). 
Moreover (see Definitions \ref{c14} and \ref{c23}(3)), we also define
$a_\delta = \bigcup\limits_{s\in [B]^{<\aleph_0}}N_s$, $e_\delta = 1$,
$\bfq' =\bfq \caret \LL \bbP^*_\delta,\Name\bbQ_\delta,a_\delta,e_\delta\RR$
and prove $\bfq' \in \gK^1$. We let $\name d(u) \defeq \name d_{|u|}(u)$.

Let $F : \omega_1 \to {}^{\omega}2$ be one-to-one such that 
$(\forall \eta \in {}^{\omega >}2)(\exists ^{\aleph_1}\alpha < \omega_1) [\eta \lhd F(\alpha)]$. (This will not be the needed 
    $\Name F$, just notation). 

For $s,t \in [B]^{<\aleph_0}$, we say $s \equiv^n_F t$ if 
$|s| = |t|$ and
$$(\forall \xi \in s)(\forall \zeta\in t)[\xi =
\bfh_{s,t}(\zeta) \Rightarrow F(\xi)\rest n = F(\zeta)\rest n].$$ 
Let $$I_n = I_n(F) \defeq \big\{s\in [B]^{<\aleph_0}:(\forall
\zeta\ne\xi\in s) [F(\zeta )\rest n\ne F(\xi)\rest n]\big\}.$$

We define $\bbR_n$ as follows: a sequence 
$\LL p_s : s \in I_n\RR\in \bbR_n$ if and only if

\begin{enumerate}
\item[(i)] for $s\in I_n$, $p_s\in \bbP^*_\lambda\cap N_s$,
\sn
\item[(ii)] for some $c_s$ we have $p_s \Vdash$ ``$\name d(s) = c_s$",
\sn
\item[(iii)] for $s,t\in I_n$, $s\equiv^n_F t\Rightarrow \bfh_{s,t}(p_t) = p_s$,
\sn
\item[(iv)] for $s,t\in I_n$, $p_s\rest N_{s\cap t} = p_t\rest N_{s\cap t}$. 
\end{enumerate}
$\bbR_n^\supminus$ is defined similarly, omitting (ii).

For $x = \LL p_s : s \in I_n\RR$ let $n(x) = n$, $p^x_s = p_s$, and (if
defined) $c^x_s = c_s$. Note that we could replace $x \in \bbR_n$ by a
finite subsequence. Let $\bbR = \bigcup\limits_{n<\omega} \bbR_n$, $\bbR^\supminus = \bigcup\limits_{n<\omega} \bbR_n^\supminus$. We define an 
order on $\bbR^\supminus$: $x \le y$ if and only if $n(x) \le n(y)$ and $$s\in I_{n(x)} \wedge t\in I_{n(y)} \wedge s \subseteq t \Rightarrow p^x_s \le p^y_t.$$

\mn
\textbf{Stage D:} Note the following facts:

\noindent\textbf{Subfact D$(\alpha)$:} If $x \in \bbR_n^\supminus$, $t\in I_n$
and $p^x_t \le p^1\in \bbP^*_\delta\cap N_t$,
\underline{then} there is $y$ such that $x \le y\in \bbR_n^\supminus$ and
$p^y_t = p^1$.

\begin{PROOF}{\bfD(\alpha)} For $s\in I_n$, we let
$$ p^y_s = \& \left\{ \bfh_{s_1,t_1}(p^1 \rest N_{t_1}) : s_1 \subseteq s,\ t_1 \subseteq t,\ s_1 \equiv^n_F t_1 \right\}\, \&\, p^x_s.$$ 

(This notation means that $p^y_s$ is a function whose domain is the union
of the domains of the conditions mentioned, and for each 
coordinate we take the canonical upper bound; see preliminaries.) 

Why is $p^y_s$ well defined? Suppose $\beta\in N_s \cap \lambda$
(for $\beta \in\lambda \setminus N_s$, clearly $p^y_s(\beta ) = \varnothing_\beta)$, $s_\ell
\subseteq s$, $t_\ell \subseteq t$, $s_\ell \equiv^n_F t_\ell$
for $\ell = 1,2$ and 
    $\beta \in \dom\big(\bfh_{s_\ell,t_\ell}(p^1\rest N_{t_\ell})\big)$, 
and it suffices to show that $p^x_s(\beta)$,
$\bfh_{s_1,t_1}(p^1\rest N_{t_1})(\beta)$, and
$\bfh_{s_2,t_2}(p^1\rest N_{t_2})(\beta)$ are pairwise
comparable. Let $u =\bigcap \big\{v \in [B]^{<\aleph_0} : \beta \in
N_v\big\}$; necessarily $u \subseteq s_1\cap s_2$, and let 
$u_\ell = \bfh^{-1}_{s_\ell,t_\ell}(u)$. As
$s_\ell,t_\ell,t\in I_n$, $s_\ell\equiv^n_F t_\ell$
and $u_\ell\subseteq t_\ell\subseteq t$,
necessarily $u_1 = u_2$. Thus $\gamma =
\bfh^{-1}_{u,v}(\beta)=\bfh^{-1}_{s_\ell,t_\ell}(\beta)$ and so the last two 
conditions are equal. 

Now $$p^x_s(\beta) = p^x_u(\beta) =\bfh_{u,v}(p^x_s(\gamma)) \le 
\bfh_{s_\ell,t_\ell}\! \big((p^x_t\rest N_{t_\ell})(\gamma)\big) =
\big(\bfh_{s_\ell,t_\ell}(p^x_t\rest N_{t_\ell})\big)(\beta).$$ 
We leave to the reader checking the other requirements.
\end{PROOF}

\sn
\textbf{Subfact D$(\beta)$:} If $x\in \bbR_n^\supminus$, $t\in I$ then
$ \bigcup\big\{p^x_s : s \in I_n,\, s \subseteq t\big\}$
(as a union of functions) exists and belongs to $\bbP^*_\lambda \cap N_t$.

\begin{PROOF}{\bfD(\beta)} See (iv) in the definition of $\bbR_n^\supminus$.
\end{PROOF}

\sn
\textbf{Subfact D$(\gamma)$:} If $x \le y$, $x\in \bbR_n$, 
$y\in \bbR_n^\supminus$, \underline{then} $y\in \bbR_n$. 

\begin{PROOF}{\bfD(\gamma)} Check it.
\end{PROOF}

\sn
\textbf{Subfact D$(\delta)$:} If $x \in \bbR_n^\supminus$, $n < m$, 
\underline{then} there is $y\in \bbR_m$ with $x \le y$. 

\begin{PROOF}{\bfD(\delta)} By subfact \textbf{D}$(\beta )$ we can find 
$x^1 = \LL p^1_t : t \in I_m\RR \in \bbR_m^\supminus$ with $x \le x^1$. Repeatedly using subfact \textbf{D}$(\alpha)$,
we can increase $x^1$ (finitely many times) to get $y \in \bbR_m$. 
\end{PROOF}

\sn
\textbf{Subfact D$(\eps)$:} If $x\in \bbR_n^\supminus$, $s,t \in I_n$,
$s\equiv^n_F t$, 
$$p^x_s \le r_1\in \bbP^*_\lambda\cap N_s, \quad p^x_t \le r_2 \in \bbP^*_\lambda \cap N_t,$$ 
$(\forall \zeta \in t)\left[F(\zeta)(n) \ne F\big(\bfh_{s,t}(\zeta) \big)(n)\right]$ 
(or just $p^x_{s_1} \rest s_1 = \bfh_{s,t}(p^x_{t_1}\rest t_1)$,
where $t_1 = \big\{ \xi \in t : F(\xi)(n) = F(\bfh_{s,t}(\xi))(n) \big\}$ and
$s_1 = \{\bfh_{s,t}(\xi) : \xi \in t_1\}$),
\underline{then} there is $y \in \bbR_{n+1}$ with $x \le y$ such that 
$r_1 = p^y_s$ and $r_2 = p^y_t$. 

\begin{PROOF}{\bfD(\eps)} Left to the reader.
\end{PROOF}

\sn
\textbf{Stage E:}\footnote{We will have $T\subset {}^{\omega >}2$ 
    from \ref{b20}(2) and then want to get a subtree with as few colors 
    as possible; we can find one isomorphic to ${}^{\omega >}2$, and there restrict ourselves to $\bigcup_n T^*_n$.
} 

We define
$T^*_k \subseteq {}^{2^k\ge} 2$ by induction on $k$
as follows:
$$\begin{aligned}T^*_0 =&\ \{\LL\ \RR,\LL 1\RR\}\cr
T^*_{k+1} =&\ T_k^* \cup \big\{\nu : 2^k <
\lh(\nu) \le 2^{k+1},\,
\nu \rest 2^k\in T^*_k, \hbox{ and}\cr 
&\ [2^k \le i < 2^{k+1} \wedge\nu(i) = 1] \Rightarrow 
i = 2^k + \big(\textstyle\sum\limits_{m<2^k}\nu(i)2^m\big)\big\}.\end{aligned}$$

We define
\begin{align*} \TrEmb(k,n) \defeq \big\{ h : &\ h \hbox{ a is
    function from $T^*_k$ into } {}^{n\ge}2\cr
    &\hbox{ such that for $\nu,\rho \in T^*_k$ we have}\cr
    &\ \eta = \nu \Leftrightarrow h(\eta)=h(\nu),\cr
    &\ \eta \lhd \nu\Leftrightarrow h(\eta) \lhd h(\nu),\cr
    &\ \lh(\eta)=\lh(\nu)\Rightarrow \lh(h(\eta)) = \lh(h(\nu)),\cr
    &\ \nu =\eta \caret \LL i\RR \Rightarrow h(\nu)\big(\lh(h(\eta))\big) = i,\cr
    &\ \lh(\eta) = k \Rightarrow \lh(h(\eta)) = n \big\}.
\end{align*}

$$
\begin{aligned}\bfT(k,n) \defeq&\ \big\{\Rang(h) : h \in 
\TrEmb(k,n)\big\},\cr
\bfT(*,n) =&\ \bigcup_k\bfT(k,n),\cr
\bfT(k,*) =&\ \bigcup_n\bfT(k,n).
\end{aligned}
$$

For $T\in \bfT(k,*)$ let $n(T)$ be the unique $n$ such that
$T\in\bfT(k,n)$ and let 
$$
\begin{aligned}
B_T =& \big\{\alpha \in B : F(\alpha)\rest n(T)
\hbox{ is a maximal member of } T\big\},\cr
\fs_T =& \big\{t \subseteq B_T:\eta\in t\wedge \nu\in t
\wedge\eta\ne\nu\Rightarrow \eta\rest n(T)\ne \nu\rest n(T) \big\},\cr 
\Theta_T =& \Big\{\LL p_s : s \in \fs_T\RR: p_s\in \bbP \cap
N_s,\big[s \subseteq t \wedge \{s,t\}\subseteq \fs_T\Rightarrow p_s = p_t\rest N_s\big]\Big\}.
\end{aligned} 
$$

Furthermore, let
$$
\begin{aligned}
\Theta_k =& \bigcup\big\{\Theta_T : T\in\bfT(k,*)\big\}\cr
\Theta=&\bigcup_k\Theta_k.
\end{aligned}
$$

For $\bar p\in \Theta$, $\bfn_{\bar p} = \bfn(\bar p)$ and $T_{\bar p}$ are defined naturally.

For $\bar p,\bar q \in \Theta$, $\bar p \le \bar q$ \underline{iff} 
$\bfn_{\bar p} \le \bfn_{\bar q}$ and 
for every $s \in \fs_{T_{\bar p}}$ we have $p_s \le q_s$.
 
\mn
\textbf{Stage F:} Let $\name g : \omega \to \omega$, 
$\name g \in N_s$, $\name g$ grows fast enough relative 
to
$\LL\sigma_n : n < \omega \RR$. We define a game \textsf{Gm}. 
A play of the game lasts $\omega$ moves: in the $n^\tthh$ move 
Player \textbf{I} chooses $\bar p^n\in \Theta_n$ and a function
$h_n$ satisfying the restrictions below, and then 
Player \textbf{II} chooses $\bar q_n \in \Theta_n$ such that 
$\bar p_n \le\bar q_n$ (so $T_{\bar p_n} = T_{\bar q_n}$).
Player \textbf{I} loses the play if at any time he has no legal move;
if he never loses, he wins. The restrictions Player \textbf{I} 
has to satisfy are: 
\begin{enumerate}
    \item[(a)] For $m < n$, $\bar q_m \le \bar p_n$, $p^n_s$ forces a value to $\name g\rest(n + 1)$.
\sn
    \item[(b)] $h_n$ is a function from $[B_{T_{\bar p_n}}]^{\le g(n)}$ to $\omega$.
\sn
    \item[(c)] $m < n\Rightarrow h_n,h_m$ are compatible.
\sn
    \item[(d)] If $m < n$, $\ell < g(m)$, and $s \in [B_{T_{\bar p_n}}]^\ell$ \underline{then} $p^n_s \Vdash \name d(s) = h_n(s)$.
\sn
    \item[(e)] Let $s_1,s_2 \in \dom(h_n)$. Then $h_n(s_1) = h_n(s_2)$ whenever $s_1,s_2$ are similar over $n$, which means:
    \begin{enumerate}
        \item[(i)] $F\big(H^\OP_{s_2,s_1}(\zeta )\big) \rest \bfn[\bar p^n] = F(\zeta ) \rest \bfn[\bar p^n]$ for $\zeta \in s_1$.
\sn
        \item[(ii)] $H^\OP_{s_2,s_1}$ preserves the relations $\essp\big(F(\zeta_1),F(\zeta_2)\big) < \essp\big(F(\zeta_3),\allowbreak F(\zeta_4)\big)$ and $F(\zeta_3)\big(\essp(F(\zeta_1),F(\zeta_2))\big) = i$ (in the interesting case $\zeta_3\ne \zeta_1$, 
        we have 
        $\zeta_2$ implies $i = 0$). 
    \end{enumerate}
\end{enumerate}

\mn
\textbf{Stage G/Claim:} Player \textbf{I} has a winning strategy in this game.

\begin{PROOF}{\bfG} 
As the game is closed, it is determined, so we assume
Player \textbf{II} has a winning strategy , and eventually we shall get
a contradiction. We define by induction on $n$, $\bar r^n$ and $\Phi^n$
such that 
\begin{enumerate}
    \item[(a)] $\bar r^n\in \bbR_n$, $\bar r^n \le \bar r^{n+1}$.
\sn 
    \item[(b)] $\Phi^n$ is a finite set of initial segments of plays of the game.
\sn
    \item[(c)] In each member of $\Phi^n$, Player \textbf{II} uses his winning strategy.
\sn
    \item[(d)] If $y$ belongs to $\Phi^n$ then it has the form $\LL\bar p^{y,\ell},h^{y,\ell},\bar q^{y,\ell} : \ell \le m(y)\RR$; let $h_y = h^{y,n_y}$ and $T_y = T_{\bar q^y,m(y)}$. Also, $T_y \subseteq ^{n\ge}2$ and $q_s^{y,\ell} \le r^n_s$ for $s \in \fs_{T_y}$. 
\sn
    \item[(e)] $\Phi_n \subseteq \Phi_{n+1}$, $\Phi_n$ is closed under taking the initial segments and the empty sequence (which too is an initial segment of a play) belongs to $\Phi_0$. 
\sn
    \item[(f)] For any $y \in\Phi_n$ and $T,h$, \underline{either} for some 
    $z \in \Phi_{n+1}$, $n_z = n_y+1$, $y = z\rest (n_y+1)$,  $T_z = T$, and 
    $h_z = h$ \underline{or} Player \textbf{I} has no legal $(n_y+1)^\tthh$ 
    move $\bar p^n,h^n$ (after $y$ was played) such that $T_{\bar p^n} = T$, $h^n = h$, and $p^n_s = r^n_s$ for $s \in \fs_T$ 
    (or always $\le$ or always $\ge$). 
\end{enumerate}

There is no problem to carry the {definition}. Now 
$\LL \bar r^n_s : n < \omega \RR$ defines a function $d^*$:
if $\eta_1,\ldots,\eta_k \in {}^m2$ are distinct then 
$d^*\big(\LL \eta_1,\ldots ,\eta_k\RR\big) = c$ \underline{iff} for
every (equivalently, `some') $\zeta_1 <\cdots < \zeta_k$ from $B$, 
$\eta_\ell \lhd F(\zeta_\ell)$ and 
$$r^k_{\{\zeta_1,\ldots,\zeta_k\}}\Vdash `` \name d_k\big(\{\zeta_1,\ldots ,\zeta_k\}\big) = c".$$

Now apply \ref{b20}(2) to this coloring and get 
$T^* \subseteq {}^{\omega >}2$ as there.
Now Player \textbf{I} could have chosen initial segments of this $T^*$ (in the $n^\tthh$ move in $\Phi_n$), and we easily get a contradiction. 
\end{PROOF}

\mn
\textbf{Stage H:} We fix a winning strategy for
Player \textbf{I} (whose existence is guaranteed by stage \textbf{G}). 

We define a forcing notion $\bbQ^*$. We have $(r,y,f)\in \bbQ^*$
\underline{iff}
\begin{enumerate}
    \item[(i)] $r \in \bbP^\cn_{a_\delta}$
\sn 
    \item[(ii)] $y = \LL \bar p^\ell,h^\ell,\bar q^\ell : \ell \le m(y)\RR$ is an initial segment of a play of \textsf{Gm} in which Player \textbf{I} uses his winning strategy. 
\sn
    \item[(iii)] $f$ is a finite function from $B$ to $\{0,1\}$ such that 
    $f^{-1}(\{1\})\in \fs_{T_y}$ (where $T_y = T_{\bar q^{m(y)}}$).
\sn
    \item[(iv)] $r = q^{y,m(y)}_{f^{-1}(\{1\})}$.
\end{enumerate}
(The order is the natural one.)

\mn
\textbf{Stage I:} If $\underline{J} \subseteq \bbP^\cn_{a_\delta}$
is dense open then $\big\{(r,y,f)\in \bbQ^* : r\in \underline{J}\big\}$
is dense in $\bbQ^*$. 

\begin{PROOF}{\bfI} By \ref{c26}(1) (by the appropriate renaming).
\end{PROOF}
 
\mn
\textbf{Stage J:} We define $\bbQ_\delta$ in $\bfV^{\bbP_\delta}$ as 
$\big\{(r,y,f)\in \bbQ^*:r \in \Name{G}_{\bbP_\delta}\big\}$,
the order is as in $\bbQ^*$.

The main point left is to prove the Knaster condition for 
the partial ordered set $\bfq^* =\bfq \caret \LL \bbP_\delta, \Name\bbQ_\delta, a_\delta, e_\delta \RR$ demanded in the definition of $\gK^1$. This will follow by
3.8(3) (after you choose meaning and renamings)
as done in stages \textbf{K} and \textbf{L} below. 

\mn
\textbf{Stage K:} So let $i < \delta$, $\cf(i)\ne \aleph_1$, and
we shall prove that $\bbP_{\delta+1}^+/ \bbP_i$ satisfies the Knaster
condition. Let $p_\alpha\in \bbP^*_{\delta+1}$ for $\alpha< \omega_1$, and we should
find $p\in \bbP_i$, $p\Vdash_{\bbP_i}$``there is an unbounded
$A\subseteq\{\alpha: p_\alpha\rest i\in \Name G_{\bbP_i}\}$ such that for
any $\alpha,\beta\in A$, $p_\alpha,p_\beta$ are compatible in
$\bbP^*_{\delta+1}/\Name G_{\bbP_i}$".

\begin{PROOF}{\bfK}
Without loss of generality:
\begin{enumerate}
    \item[(a)] $p_\alpha \in \bbP^\cn_{\delta+1}$
\sn
    \item[(b)] For some $\LL i_\alpha : \alpha < \omega_1\RR$ increasing continuous with limit $\delta$ we have $i_0 > i$, $\cf(i_\alpha) \ne \aleph_1$, $p_\alpha\rest\delta \in \bbP_{i_{\alpha+1}}$, and 
    $p_\alpha\rest i_\alpha\in \bbP_{i_0}$. 
\sn
    Let $p^0_\alpha = p^\alpha\rest i_0$, $p_\alpha^1 = 
    p_\alpha \rest \delta = p_\alpha \rest i_{\alpha+1}$, and 
    $p_\alpha(\delta) = (r_\alpha, y_\alpha, f_\alpha)$.
\sn
    \item[(c)] $r_\alpha\in \bbP_{i_{\alpha+1}}$, 
$r_\alpha\rest i_\alpha\in \bbP_{i_0}$, and $m(y_\alpha)=m^*$.
\sn 
    \item[(d)] $\dom(f_\alpha) \subseteq i_0\cup [i_\alpha, i_{\alpha+1})$,
\sn 
    \item[(e)] $f_\alpha\rest i_0$ is constant. (Remember, $\otp(B)=\omega_1$.)
\sn 
    \item[(f)] If $\dom (f_\alpha) = \{j^\alpha_0,\ldots 
    j^\alpha_{k_\alpha-1}\}$ then $k_\alpha = k$, 
    $[j^\alpha_\ell < i_\alpha \Leftrightarrow \ell < k^*]$, $\bigwedge\limits_{\ell<k^*} j^\alpha_\ell=j^\ell$, $f(j_\ell^\alpha)=f(j^\beta_\ell)$, and
    $F(j_\ell^\alpha)) \rest m(y_\alpha) = F(j^\beta_\ell)\rest m(y_\beta)$. 
\end{enumerate}
The main problem is the compatibility of the $q^{y_\alpha,m(y_\alpha)}$. 
Now by the definition of $\Theta_\alpha$ (in stage \textbf{E}) and \ref{c26}(3) this holds.
\end{PROOF}

\noindent\textbf{Stage L:} If $c \subset \delta+1$ is closed for $\bfq^*$,
then $\bbP^*_{\delta+1}/ \bbP^\cn_c$ satisfies the Knaster condition.

If $c$ is bounded in $\delta$, choose a successor $i\in (\sup c,\delta)$
for $\bfq\rest i\in \gK_1$. We know that $\bbP_i/ \bbP^\cn_c$ satisfies
the Knaster condition and by stage K, $\bbP_{\delta+1}^*/ \bbP_i$ also
satisfies the Knaster condition; as it is preserved by
composition we have finished the stage.

So assume $c$ is unbounded in $\delta$ and it is easy too. So as
seen in stage \textbf{J}, we have finished the proof of \ref{c2}. 
\end{PROOF}

\begin{theorem}\label{c38} If $\lambda \ge \beth_\omega$ and
$\bbP$ is the forcing notion which adds $\lambda$ Cohen reals, then:
\begin{enumerate}
    \item[$(*)_1$] In $\bfV^\bbP$, if $n < \omega$ and 
    $d : [\lambda ]^{\le n} \to \sigma$ with $\sigma < \aleph_0$, 
    \underline{then} for some c.c.c. forcing notion $\bbQ$ we have 
    $\Vdash_\bbQ$ ``\emph{there are an uncountable $A \subseteq \lambda$ 
    and a one-to-one $F : A \to {}^\omega 2$ such that $d$ is $F$-canonical on $A$}'' (see notation in \S2).
\sn
    \item[$(*)_2$] If $\lambda \ge \mu \to_\wsp(\kappa)_{\aleph_0}$ in $\bfV$ (see \cite{Sh:289}) and $d : [\mu]^{\le n} \to \sigma$ in $\bfV^\bbP$ (with $\sigma < \aleph_0$) \underline{then}, in $\bfV^\bbP$, for some c.c.c. forcing notion $\bbQ$ we have $\Vdash_\bbQ$ \emph{``there are 
    $A \in [\mu]^\kappa$ and one-to-one $F:A \to {}^\omega 2$ such that $d$ is $F$-canonical on $A$"} (see \S2). 
\sn
    \item[$(*)_3$] If $\lambda \ge \mu \to_\wsp(\aleph_1)^n_{\aleph_2}$ in $\bfV$ and $d : [\mu]^{\le n} \to \sigma$ in $\bfV^\bbP$ (with 
    $\sigma < \aleph_0$) \underline{then}, in $\bfV^\bbP$, for every 
    $\alpha < \omega_1$ and $F : \alpha \to {}^\omega 2$, for some 
    $A \subseteq \mu$ of order type $\alpha$ and $F' : A \to {}^\omega 2$, $F'(\beta) = F(\otp(A \cap \beta))$, $d$ is $F'$-canonical on $A$. 
\sn
    \item[$(*)_4$] In $\bfV^\bbP$, $2^{\aleph_0} \to (\alpha,n)^3$ for every $\alpha < \omega_1$ and $n < \omega$. Really, assuming\\  $\bfV \models \GCH$ we have $\aleph_{n^1_3} \to (\alpha,n)$ (see \cite{Sh:289}).
\end{enumerate}
\end{theorem}

\begin{PROOF}{\ref{c38}}
Similar to the proof of \ref{c2}. Superficially we need more
indiscernibility then we get, but getting $\LL M_u:u\in
[B]^{\le n}\RR$ we ignore $d(\{\alpha ,\beta\})$ when there is no $u$
with $\{\alpha ,\beta\} \in M_u$.
\end{PROOF}

\begin{theorem}\label{c41} 
If $\lambda$ is strongly inaccessible $\omega$-Mahlo and $\mu < \lambda$, 
then for some c.c.c. forcing notion $\bbP$ of cardinality $\lambda$, 
$\bfV^\bbP$ satisfies
\begin{enumerate}[$(a)$]
    \item $\MA_\mu$
\sn 
    \item $2^{\aleph_0} = \lambda = 2^\kappa$ for $\kappa < \lambda$.
\sn 
    \item $\lambda \to [\aleph_1]^n_{\sigma,h(n)}$ for $n < \omega$, 
    $\sigma < \aleph_0$, and $h(n)$ as in \ref{c2}.
\end{enumerate}
\end{theorem}

\begin{PROOF}{\ref{c41}} 
Again, like \ref{c2}.
\end{PROOF}

\newpage

\section{Partition theorem for trees on large cardinals}

\begin{lemma}\label{d2} 
Suppose $\mu > \sigma + \aleph_0$ and
\begin{enumerate}
    \item[$(*)_\mu$] for every $\mu$-complete forcing notion $\bbP$, in
$\bfV^\bbP$, $\mu$ is measurable.  
\end{enumerate}

\noindent
Then
\begin{enumerate}[$(1)$]
    \item We have $\Pr^\fe_\eht(\mu ,n,\sigma)$ for all $n < \omega$.
\sn
    \item $\Pr^\fe_\eht(\mu ,<\aleph_0,\sigma)$, if there is 
    $\lambda > \mu$ such that $\lambda \to \big(\mu^+\big)^{<\omega}_2$. 
\sn
    \item In both cases we can have the $\Pr^\fe_\ehtn$ version, and even choose the\\ $\LL {<^*_\alpha} : \alpha < \mu \RR$ in any of the following ways. 
    \begin{enumerate}[$(a)$]
        \item We are given $\LL <^0_\alpha : \alpha < \mu \RR$, and (for $\eta,\nu \in {}^\alpha 2\cap T$, $\alpha\in \SP(T)$, and $T$ the subtree we consider) we let: 
        \begin{itemize}
            \item $\eta <^*_\alpha \nu$ if and only if 
            $\clp_T(\eta) <^0_\beta \clp_T(\nu)$, where 
            $\beta =$\\ $\otp(\alpha \cap \SP(T))$ and 
            $\clp_T(\eta) = \big\LL \eta(j) : j \in \lh(\eta),\ j \in \SP(T) \big\RR$.
        \end{itemize}
\sn
        \item We are given $\LL <^0_\alpha : \alpha < \mu \RR$, and we
        say $\eta <^*_\alpha\nu$ if and only if\\ $n \rest (\beta +1) <^0_{\beta+1} \nu \rest (\beta+1)$, where 
        $\beta = \sup(\alpha \cap \SP(T))$. 
    \end{enumerate}
\end{enumerate}
\end{lemma}

\begin{remark}\label{d5}
{1) $(*)_\mu$ holds for a supercompact after Laver
treatment. On hypermeasurable, see Gitik-Shelah \cite{Sh:344}.

\sn
2) We can in $(*)_\mu$ restrict ourselves to the forcing notion $\bbP$
actually used. For that, by Gitik \cite{Gi10}
much smaller large cardinals suffice.

\sn
3) The proof of \ref{d2} is a generalization of a proof of
Harrington to the  Halpern-Lauchli theorem from 1978.}
\end{remark}

\begin{conclusion}\label{d8} In \ref{d2} we can get 
$\Pr^\fe_\hgt(\mu ,n,\sigma)$ (even with (3)).
\end{conclusion}

\begin{PROOF}{\ref{d8}} 
We do the parallel to \ref{d2}(1).
By $(*)_\mu$, $\mu$ is weakly compact
hence by \ref{b17}(2) it is enough to prove 
$\Pr^\fe_\aht(\mu,n,\sigma)$. This
follows from \ref{d2}(1) by \ref{b17}(1). 
\end{PROOF}

\begin{PROOF}{\ref{d2}} \textbf{Proof of \ref{d2}:} 

\sn
1), 2). Let $\kappa \le \omega$,
$\sigma(n) < \mu$, $d_n \in \Col^n_{\sigma(n)}(^{\mu >}2)$ for 
$n < \kappa$.

Choose $\lambda$ such that $\lambda \to (\mu^+)^{<2\kappa}_{2^\mu}$
(there is such a $\lambda$ by assumption for (2) and 
by $\kappa < \omega$ for (1)). Let $\bbQ$ be the forcing notion
$(^{\mu >}2,\lhd)$, and $\bbP = \bbP_\lambda$ be 
$$\big\{ f : \dom(f) \text{ is a subset of $\lambda$ of cardinality} <\mu,\ f(i) \in \bbQ \big\},$$ 
ordered naturally.
For $i \notin \dom(f)$, take $f(i)= \LL\ \RR$. 
Let $\name \eta_i$ be the $\bbP$-name for $\bigcup \{f(i) : f \in \Name G_\bbP\}$. Let 
$\Name D$ be a $\bbP$-name of a normal ultrafilter over $\mu$. 
For each $n < \omega$, $d \in \Col^n_{\sigma(n)}(^{\mu >}2)$,
$j < \sigma(n)$ and
$u= \{\alpha_0,\ldots ,\alpha_{n-1}\}$, where $\alpha_0 <\cdots <
\alpha_{n-1} < \lambda$, let 
$\Name{A}^j_d(u)$ be the $\bbP_\lambda$-name of the set 
$$A^j_d(u) = \Big\{i < \mu :\LL
\name\eta_{\alpha_\ell}\rest i :\ell < n \RR \text{ are pairwise 
distinct, } j = d(\eta_{\alpha_0}\rest i,\ldots,\eta_{\alpha_{n-1}}\!\rest i)\Big\}.$$
So $\Name A^j_d(u)$ is a $\bbP_\lambda$-name of a subset of
$\mu$, and for $j(1) < j(2) < \sigma(n)$ we have 
$\Vdash_{\bbP_\lambda}$``$\Name A^{j(1)}_d(u)\cap \Name A^{j(2)}_d(u) = \varnothing$, and
$\bigcup_{j<\sigma(n)} \Name A^j_d(u)$ is a co-bounded subset of
$\mu$''. As $\Vdash_\bbP$ ``$\gD$ is $\mu$-complete uniform
ultrafilter on $\mu$", in $\bfV^\bbP$ there is exactly one 
 $j < \sigma(n)$ with $ A^j_d(u)\in\gD$. Let
$\name j_d(u)$ be the
$\bbP$-name of this $j$. 

Let $I_d(u) \subseteq \bbP$ be a maximal antichain of $\bbP$, each
member of $I_d(u)$ forces a value to $\name j_d(u)$. Let
$W_d(u) = \bigcup \{\dom(p) : p \in I_d(u)\}$ and
$W(u) =\bigcup \{W_{d_n}(u) : n < \kappa \}$. So 
$W_d(u)$ is a subset of $\lambda$ of cardinality $\le \mu$ as
well as $W(u)$ (as $\bbP$ satisfies the $\mu^+$-c.c. and $p\in
P\Rightarrow |\dom(p)| < \mu )$. 

As $\lambda \to (\mu^{++})^{<2\kappa}_{2^\mu}$, 
$d_n \in \Col^n_{\sigma_n}(^{\mu >}2)$ there is a subset $Z$ of $\lambda$ of
cardinality $\mu^{++}$ and set $W^+(u)$ for each
$u \in [Z]^{<\kappa}$ such that: 
\begin{enumerate}
    \item[(i)] $W^+(u_1)\cap W^+(u_2) = W^+(u_1\cap u_2)$
\sn
    \item[(ii)] $W(u) \subseteq W^+(u)$ if $u\in [Z]^{<\kappa}$.
\sn
    \item[(iii)] If $|u_1| = |u_2| < \kappa$ and $u_1,u_2 \subseteq Z$ then 
    $W^+(u_1)$ and $W^+(u_2)$ have the same order type. 
    
    (Note that $H[u_1,u_2] = H^\OP_{W^+(u_1),W^+(u_2)}$ naturally induces a map from $\bbP\rest u_1 = \{p\in \bbP : \dom(p) \subseteq W^+(u_1)\}$ to\\ $\bbP\rest u_2 = \{p \in \bbP : \dom(p) \subseteq W^+(u_2)\}$.)
\sn
    \item[(iv)] If $u_1,u_2 \in [Z]^{<\kappa}$ and $|u_1| = |u_2|$ 
    then $H[u_1,u_2]$ maps $I_{d_n}(u_1)$ onto $I_{d_n}(u_2)$ and 
    $$q \Vdash ``{\name j}_d(u_1) = j" \Leftrightarrow H[u_1,u_2](q) 
    \Vdash ``{\name j}_d(u_2) = j".$$ 
    \item[(v)] If $u_1 \subseteq u_2\in [Z]^{<\kappa}$, 
    $u_3 \subseteq u_4 \in [Z]^{<\kappa}$, $|u_4| = |u_2|$, and
    $H^\OP_{u_2,u_4}$ maps $u_1$ onto $u_3$, \underline{then}
    $H[u_1,u_3] \subseteq H[u_2,u_4]$. 
\end{enumerate}

Let $\gamma(i)$ be the $i^\tthh$ member of $Z$.

Let $s(m)$ be the set of the first $m$ members of $Z$ and 
$$
\bbR_n = \big\{p\in \bbP:\dom(p) \subseteq W^+(s(n))\ 
{\setminus} \bigcup_{t\subset
s(n)}W^+(t)\big\}.
$$ 
We define, by induction on $\alpha < \mu$, a function $F_\alpha$ and 
$p_u \in \bbR_{|u|}$ for $u \in \bigcup\limits_{\beta <\alpha}[^\beta2]^{<\kappa}$ 
where we let $\varnothing_\beta$ be the empty subset of $[^\beta 2]$, 
we behave as if $[\beta \ne \gamma \Rightarrow \varnothing_\beta\ne \varnothing_\gamma]$, and we also define
$\zeta(\beta) < \mu$ such that: 
\begin{enumerate}
    \item[(i)] $F_\alpha$ is a function from ${}^{\alpha >}2$ into ${}^{\mu >}2$,
    extending $F_\beta$ for each $\beta < \alpha$.
\sn
    \item[(ii)] $F_\alpha$ maps ${}^\beta 2$ to ${}^{\zeta(\beta)}2$ for some $\zeta(\beta) < \mu$, and $$\beta_1 < \beta_2 < \alpha \Rightarrow \zeta (\beta_1) < \zeta (\beta_2).$$
    \item[(iii)] $\eta \lhd \nu \in {}^{\alpha >}2$ implies $F_\alpha(\eta) \lhd F_\alpha (\nu)$.
\sn
    \item[(iv)] For $\eta\in ^{\beta}2$, $\beta + 1 < \alpha$, and $\ell< 2$, we have $F_\alpha(\eta) \caret \LL\ell\RR \unlhd F_\alpha(\eta \caret \LL\ell\RR)$. 
\sn
    \item[(v)] $p_u\in \bbR_m$ whenever $u\in[^\beta 2]^m$, $m < \kappa$,
    $\beta < \alpha$ and for $u(1)\in [Z]^m$ let $p_{u,u(1)} = H[s(|u|),u(1)](p_u)$. 
\sn
    \item[(vi)] $\eta \in {}^\beta 2$, $\beta < \alpha$, then 
    $p_{\{\eta\}}(\min Z) = F_\alpha(\eta)$.
\sn
    \item[(vii)] If $\beta < \alpha$, $u \in [^\beta2]^n$, $n < \kappa$, and
    $h : u \to s(n)$ is one-to-one and onto (but not necessarily order preserving) \underline{then} for some $c(u,h) < \sigma(n)$,
    $$
    \bigcup_{t \subseteq u}p_{t,h''(t)} \Vdash_{\bbP_\lambda} ``\name d{}_n( \hbox{$\name\eta_{\gamma(0)}$},\ldots,\hbox{$\name\eta_{\gamma(n-1)}$}) = c(u,h)".
    $$
    (Note: as $p_u \in \bbR_{|u|}$, the domains of the conditions in this union are pairwise disjoint.)
\sn
    \item[(viii)] If $n,u,\beta ,h$ are as in (vii), 
    $u = \{\nu_0,\ldots ,\nu_{n-1}\}$, $\nu_\ell \lhd \rho_\ell \in {}^\gamma 2$, 
    and $\beta \le \gamma < \alpha$, then $d_n(F_\alpha (\rho_0),\ldots , F_\alpha (\rho_{n-1})) = c(u,h)$, where $h$ is the unique function from $u$ onto $s(n)$ such that 
    $[h(\nu_\ell) \le h(\nu_m) \Rightarrow \rho_\ell <^*_\gamma \rho_m]$.
\sn
    \item[(ix)] If $\beta < \gamma < \alpha$, $\nu_1,\ldots, \nu_{n-1} \in {}^\gamma 2$, $n <\kappa$, and $\nu_0\rest\beta, \ldots, \nu_{n-1}\rest\beta$ are pairwise distinct, then: 
    $p_{\{\nu_0\rest\,\beta,\ldots ,\nu_n\rest\,\beta \}} \subseteq 
    p_{\{\nu_0,\ldots ,\nu_{n-1}\}}$.
\end{enumerate}

\sn
\textbf{For $\alpha$ limit}: no problem.

\sn
\textbf{For $\alpha+1$ with $\alpha$ limit}: we try to define
$F_\alpha(\eta)$ for $\eta\in {}^\alpha 2$ such that 
$$\bigcup_{\beta <\alpha} F_{\beta +1}(\eta\rest\beta) \unlhd F_\alpha (\eta)$$ 
and (viii) holds. Let $\zeta = \bigcup\limits_{\beta < \alpha}\zeta(\beta)$.
For $\eta \in {}^\alpha 2$, we define 
$$F^0_\alpha(\eta) \defeq \bigcup_{\beta <\alpha} F_\alpha(\eta\rest\beta)$$
and for $u\in[^\alpha 2]^{<\kappa}$, 
$$p^0_u = \bigcup \big\{p^0_{\{\nu\rest\beta : \nu \in u\}} : \beta < \alpha \wedge \big| \{\nu\rest\beta:\nu\in u\}\big| = |u|\big\}.$$
Clearly $p^0_u\in \bbR_{|u|}$.

Then let $h : {}^\alpha 2 \to Z$ be one-to-one such that
$\eta <^*_\alpha \nu \Leftrightarrow h(\eta) < h(\nu)$ and let 
$$p = \bigcup\{p^0_{u,u(1)} : u(1) \in [Z]^{<\kappa},\ 
u \in [{}^\alpha 2]^{<\kappa},\ |u(1)| = |u|,\ h''(u) = u(1)\}.$$

For any generic $G \subseteq \bbP_\lambda$ to which $p$ belongs, for 
$\beta < \alpha$, $n < \omega$, and ordinals $i_0 <\cdots <i_{n-1}$ 
from $Z$ such that $\LL h^{-1}(i_\ell)\rest \beta : \ell < n\RR$ 
are pairwise distinct, we have that 
$$B_{\{i_\ell : \ell < n\},\beta} \defeq \Big\{ \xi < \mu : d_n(\eta_{i_0}\rest\xi,\ldots,\eta_{i_{n-1}}\rest\xi) = c(u,h^*) \Big\}$$
belongs to $\gD[G]$, where 
$u = \{h^{-1}(i_\ell)\rest\beta : \ell < n\}$ and
$ h^* : u \to s(|u|)$ is defined by 
$h^*(h^{-1}(i_\ell)\rest\beta) =
H^\OP_{\{i_\ell : \ell<n\},s(n)}(i_\ell)$. 
Really every large enough $\beta < \mu$
can serve so we omit it. As $\gD[G]$ is $\mu$-complete uniform
ultrafilter on $\mu$, we can find $\xi \in (\zeta,\kappa)$ such
that $\xi \in B_u$ for every $u \in [^\alpha 2]^{<\omega}$. 

For $\nu \in {}^\alpha 2$, we let 
$F_\alpha (\nu) = \name \eta_{h(i)}[G]\rest\xi$, and we let 
$p_u = p^0_u$ except when $u = \{\nu\}$. In that case: 
$$ p_u(i) = 
\begin{cases}
    p^0_u(i)&          \text{if } i \ne \gamma(0)\cr
    F_{\alpha+1}(\nu)& \text{if } i =   \gamma(0).
\end{cases}$$

\sn
\textbf{For $\alpha +1$, with $\alpha$ a successor}: 

First, for $\eta\in {}^{\alpha-1}2$ define
$F\big(\eta \caret \LL\ell\RR\big) = F_\alpha (\eta) \caret \LL\ell\RR$.
Next we let $\big\{(u_i,h_i):i < i^*\big\}$ list all
pairs $(u,h)$ with $u\in [^\alpha 2]^{\le n}$ and 
$h : u \to s(|u|)$ one-to-one and onto. Now, by
induction on $i \le i^*$, we define 
$p^i_u$ (for $u\in [^\alpha 2]^{<\kappa}$) such that:
\begin{enumerate}
    \item[(a)] $p^i_u\in \bbR_{|u|}$
\sn
    \item[(b)] $p^i_u$ increases with $i$.
\sn
    \item[(c)] For $i + 1$, clause (vii) above holds (with $\alpha,u_i,h_i$ here standing in for $\beta,u,h$ there). 
\sn
    \item[(d)] If $\nu_m \in {}^\alpha 2$ for $ m < n < \kappa$ and $\LL\nu_m\rest(\alpha-1) : m < n\RR$ are pairwise distinct, \underline{then} $p_{\{\nu_m \rest (\alpha-1) : m < n\}} \le p^0_{\{\nu_m : m < n\}}$.
\sn
    \item[(e)] If $\nu \in {}^\alpha2$ and $\nu(\alpha-1) = \ell$ 
    \underline{then} $p^0_{\{\nu\}}(0) = F_\alpha(\nu\rest(\alpha - 1)) \caret \LL\ell\RR$. 
\end{enumerate}

There is no problem to carry the induction.

Now $F_{\alpha +1}\rest {}^\alpha 2$ is to be defined as in the second
case, starting with $\eta \to p^{i^*}_{\{\eta\}}(\eta)$. 
 
\sn
\textbf{For $\alpha = 0,1$:} Left to the reader.

So we have finished the induction
hence the proof of $4.1(1),(2)$.

\sn
3) Left to the reader (the only influence is the choice
of $h$ in stage of the induction).
\end{PROOF}

\newpage

\section{Somewhat complementary negative partition relation in ZFC}

The negative results here suffice to show that the value we have
for $2^{\aleph_0}$ in \S3 is reasonable. In particular, the Galvin
conjecture is wrong and that for every $n < \omega$, for some 
$m < \omega$, $\aleph_n \not\to [\aleph_1]^m_{\aleph_0}$. 

See Erd\H{o}s-Hajnal-M{\'a}t{\'e}-Rado \cite{EHMR} for

\begin{fact}\label{e2} 
If $2^{<\mu} < \lambda \le 2^\mu$ and $\mu \not\to [\mu]^n_\sigma$ then 
$\lambda \not\to [(2^{<\mu})^+]^{n+1}_\sigma$.
\end{fact}

This shows that if e.g. in \ref{a11} we want to increase the
exponents to 3 (and still $\mu = \mu ^{<\mu}$) {then} $\mu$ cannot
be successor (when $\sigma \le \aleph_0$; by \cite[3.5(2)]{Sh:276}). 

\begin{definition}\label{e5} $\Pr_\np(\lambda ,\mu ,\olsi\sigma)$ 
(where $\olsi\sigma = \LL\sigma_n : n < \omega \RR$) means that there are 
functions $F_n : [\lambda]^n \to \sigma_n$ such that for every 
$W \in [\lambda]^\mu$, for some $n$, $F''_n([W]^n) = \sigma(n)$. 
The negation of this property is denoted by
$\NPr_\np(\lambda ,\mu,\olsi\sigma)$.
\end{definition}

If the sequence in constantly $\sigma$ we may write $\sigma$ instead of 
$\LL\sigma_n : n < \omega \RR$.

\begin{remark}\label{e8} 
1) Note that $\lambda \to [\mu ]^{<\omega}_\sigma$ means ``if
$F:[\lambda]^{<\omega} \to \sigma$ then for some $A \in [\lambda ]^\mu$,
$F''([A]^{<\omega}) \ne \sigma$." So for $\lambda \ge \mu \ge \sigma = \aleph_0$, {we have} $\lambda \not\to [\mu]^{<\omega}_\sigma$ 
(use $\alpha \mapsto |\alpha|$ for $F$), and
$\Pr_\np(\lambda ,\mu ,\sigma)$ is stronger than 
$\lambda \not\to [\mu ]^{<\omega}_\sigma$.

\sn
2) We do not write down the monotonicity properties of $\Pr_\np$:
they are obvious.
\end{remark}

\begin{claim}\label{e11} 
1) Without loss of generality we can (in \ref{e5}) use 
$F_{n,m} : [\lambda ]^n \to \sigma_m$
for $n,m <\omega$ and obvious monotonicity properties holds, and
$\lambda \ge \mu \ge n$.

\sn
2) Suppose $\NPr_\np(\lambda ,\mu ,\kappa)$ and $\kappa\not\to
[\kappa]^n_\sigma$, or even $\kappa \not\to [\kappa]^{<\omega}_\sigma$. 
\underline{Then} the following case of the Chang conjecture holds: 
\begin{enumerate}
    \item[$(*)$] For every model $M$ with universe $\lambda$ and countable vocabulary, there is an elementary submodel $N$ of $M$ of cardinality $\mu$ 
    with
\end{enumerate}

\noindent
3) If $\NPr_\np(\lambda ,\aleph_1,\aleph_0)$ then
$(\lambda,\aleph_1) \to (\aleph_1,\aleph_0)$.
\end{claim}

\begin{PROOF}{\ref{e11}} 
Easy.
\end{PROOF}

\begin{theorem}\label{e14} 
Suppose $\Pr_\np(\lambda_0,\mu,\aleph_0)$, $\mu$ is regular $ >\aleph_0$ 
and $\lambda_1 \ge \lambda_0$, and no 
$\mu' \in (\lambda_0,\lambda_1)$ is $\mu'$-Mahlo. Then 
$\Pr_\np(\lambda_1,\mu,\aleph_0)$.
\end{theorem}

\begin{PROOF}{\ref{e14}} 
Let $\chi = \beth_8 (\lambda_1)^+$, let 
$\{F^0_{n,m} : m < \omega\}$ list the definable $n$-place 
functions in the model $(\clH(\chi),\in,<^*_\chi)$ 
with $\lambda_0,\mu ,\lambda_1$ as parameters, let 
$F^1_{n,m}(\alpha_0, \ldots, \alpha_{n-1})$
(for $\alpha_0,\ldots ,\alpha_{n-1} < \lambda_1$) be equal to 
$F^0_{n,m}(\alpha_0,\ldots ,\alpha_{n-1})$ if it is an ordinal 
$< \lambda_1$ and zero otherwise. Let $F_{n,m}(\alpha_0,\ldots,\alpha_{n-1})$ 
(for $\alpha_0,\ldots ,\alpha_{n-1} < \lambda_1$)
be $F^0_{n,m}(\alpha_0,\ldots ,\alpha_{n-1}$) if it is an ordinal $<\omega$ and
zero otherwise. We shall show that the $F_{n,m}$ (for $n,m < \omega$) exemplify
$\Pr_\np(\lambda_1,\mu ,\aleph_0)$ (see \ref{e8}(1)). 

So suppose $W\in [\lambda_1 ]^\mu$ is a counterexample to
$\Pr(\lambda_1,\mu ,\aleph_0)$: i.e. for no $n,m$ is 
$F''_{n,m}([W]^n) = \omega$. 
Let $W^*$ be the closure of $W$ under $\{F^1_{n,m} : n,m < \omega \}$. 
Let $N$ be the Skolem Hull of $W$ in $(\clH(\chi),\in,<^*_\chi)$, 
so clearly $N \cap \lambda_1 = W^*$. 
(Note $W^* \subseteq \lambda_1$ 
and
$|W^*| = \mu$.) 
Also, as $\cf(\mu) > \aleph_0$, if $A \subseteq W^*$ with $|A| = \mu$
then for some $n,m < \omega$ and $u_i \in [W]^n$ (for $i < \mu$)
we have $F^1_{n,m}(u_i)\in A$ and 
$$i < j < \mu \Rightarrow F^1_{n,m}(u_i) \ne F^1_{n,m}(u_{
{j}}).$$ 
It is easy to check that also $W^1 \defeq \{F^1_{n,m}(u_i) : i < \mu \}$ 
is a counterexample to $\Pr(\lambda_1,\mu,\sigma)$. In particular, for 
$n,m < \omega$, $W_{n,m} = \{F^1_{n,m}(u):u\in [W]^n\}$ is a 
counterexample if it has power $\mu$. Without loss of generality $W$
is a counterexample with minimal 
$\delta \defeq \sup(W) = \bigcup\{\alpha + 1 : \alpha \in W\}$.
The above discussion shows that $|W^* \cap \alpha| < \mu$ for $\alpha < \delta$. Obviously $\cf(\delta) = \mu^+$. Let $\LL \alpha_i : i < \mu \RR$ be a strictly 
increasing sequence of members of $W^*$, converging to $\delta$, 
such that for limit $i$ we have 
$\alpha_i = \min \! \big( W^* \setminus \bigcup\limits_{j<i}(\alpha_j + 1) \big)$. 
Let $N = \bigcup\limits_{i < \mu} N_i$ where $N_i\prec N$, $|N_i| < \mu$, $N_i$ 
increasing continuous, and without loss of generality 
$N_i \cap \delta = N \cap \alpha_i$. 

\sn
\textbf{Fact $(\alpha)$}: $\delta > \lambda_0$.

\begin{PROOF}{\alpha}
Otherwise we then get an easy contradiction to $\Pr(\lambda_0,\mu,\sigma)$, 
as when choosing the $F^0_{n,m}$ we allowed $\lambda_0$ as a parameter. 
\end{PROOF}
 
\sn
\textbf{Fact $(\beta)$}: If $F$ is a unary function
definable in $N$, $F(\alpha)$ is a club of $\alpha$ for every limit
ordinal $\alpha\ (<\lambda_1)$ \underline{then} for some club $C$ of 
$\mu$ we have 
$$\big(\forall j \in C \setminus \{\min C\}\big)
\big(\exists i_1< j\big)\big(\forall i\in (i_1,j)\big)
\big[i\in C\Rightarrow\alpha_i \in F(\alpha_j)\big].$$ 

\begin{PROOF}{\beta}
For some club $C_0$ of $\mu$ we have 
$$j\in C_0 \Rightarrow \big( N_j, \{ \alpha_i : i < j\}, W \big) \prec 
\big(N, \{\alpha_i : i < \mu\}, W \big).$$
We let $C = C'_0 = \acc(C)$ ($=$ set of accumulation points of $C_0$).

We check $C$ is as required; suppose $j$ is a counterexample.
So $j = \sup(j\cap C)$ (otherwise choose $i_1 = \max(j\cap C)$). So we
can define, by induction on $n$, a sequence of $i_n$ such that: 
\begin{enumerate}
    \item[(a)] $i_n < i_{n+1} < j$
\sn 
    \item[(b)] $\alpha_{i_n} \notin F(\alpha_j)$
\sn 
    \item[(c)] $(\alpha_{i_n},\alpha_{i_{n+1}}) \cap F(\alpha_j) \ne \varnothing $.
\end{enumerate} 

Why $(C'_0)$? $\models$ ``$F(\alpha_j)$ is unbounded below $\alpha_j$" hence
$N\models$ ``$F(\alpha_j)$ is unbounded below $\alpha_j$", but in $N$,
$\{\alpha_i : i\in C_0,\,\,i < j\}$ is unbounded below $\alpha_j$. 

Clearly, for some $n,m$ we have $\alpha_j\in W_{n,m}$ (see above). Now we
can repeat the proof of \cite[3.3(2)]{Sh:276}\footnote{See mainly the end.} using
only members of $W_{n,m}$.\\
Note: here we set the number of colors to be $\aleph_0$.
\end{PROOF}
 
\sn
\textbf{Fact $(\beta)^+$}: Without loss of generality, the club $C$ in Fact $(\beta)$ is $\mu$.

\sn
Proof: By renaming.

\mn
\textbf{Fact $(\gamma)$}: $\delta$ is a limit cardinal.

\begin{PROOF}{\gamma}
Suppose not. Now $\delta$ cannot be a successor cardinal 
(as $\cf(\delta) = \mu \le \lambda_0 <\delta$) hence for every large enough $i$,
$|\alpha_i| = |\delta|$, so $|\delta| \in W^* \subseteq N$ and
$|\delta|^+\in W^*$. 

So $W^*\cap |\delta|$ has cardinality $<\mu$ hence order-type
equal to some $\gamma^* < \mu$. Choose $i^* < \mu$ limit such that 
$[j < i^* \Rightarrow j + \gamma^* < i^*]$. There is a definable 
function $F$ of\\ $(\clH(\chi),\in,<^*_\chi)$ such that for every limit 
ordinal $\alpha$, $F(\alpha)$ is a club of $\alpha$, such that if 
$|\alpha| < \alpha$ then $F(\alpha) \cap |\alpha| = \varnothing$ and
$\otp(F(\alpha)) = \cf(\alpha)$. 

So in $N$ there is a closed unbounded subset $C_{\alpha_j} =
F(\alpha_j)$ of $\alpha_j$ of order type $\le \cf(\alpha_j) \le |\delta|$, hence
$C_{\alpha_j}\cap N$ has order type $\le\gamma^*$, hence for $i^*$
chosen above unboundedly many $i < i^*$, $\alpha_i\notin
C_{\alpha_{i^*}}$. We can finish by Fact $(\beta)^+$. 
\end{PROOF} 
 
\sn
\textbf{Fact $(\delta)$}: For each $i < \mu$, $\alpha_i$ is a
cardinal. 

\begin{PROOF}{\delta}
If $|\alpha_i|<i$ then $|\alpha_i|\in N_i$, but then
$|\alpha_i|^+\in N_i$ contradicting Fact ($\gamma$), by which
$|\alpha_i|^+ <\delta$, as we have assumed $N_i\cap\delta= N\cap\alpha_i$. 
\end{PROOF} 
 
\sn
\textbf{Fact $(\eps)$}: For a club of $i < \mu$,
$\alpha_i$ is a regular cardinal.

\begin{PROOF}{\eps}
If $S = \{i : \alpha_i$ singular$\}$ is stationary, then
the function $\alpha_i \mapsto \cf(\alpha_i)$ is regressive on $S$. 
By Fodor's lemma, for some $\alpha^* <\delta$, 
$\{i < \mu : \cf (\alpha_i) < \alpha^*\}$ is
stationary. As $|N\cap\alpha^*| < \mu$ for some $\beta^*$, 
$\{i < \mu : \cf (\alpha_i) = \beta^*\}$ is stationary. Let $F_{1,m}(\alpha)$ be a
club of $\alpha$ of order type $\cf(\alpha)$, and by Fact ($\beta$) we get
a contradiction as in Fact ($\gamma$). 
\end{PROOF} 

\sn
\textbf{Fact $(\zeta)$}: For a club of $i < \mu$, $\alpha_i$ is Mahlo.

\begin{PROOF}{\zeta}
Use $F_{1,m}(\alpha) =$ a club of $\alpha$ which, if $\alpha$ is
 a successor cardinal or 
inaccessible not Mahlo, then it contains no inaccessible, and
continue as in Fact ($\gamma$). 
\end{PROOF}
 
\sn
\textbf{Fact $(\xi)$}: For a club of $i < \mu$, $\alpha_i$
is $\alpha_i$-Mahlo.

\begin{PROOF}{\xi}
Let $F_{1,m(0)}(\alpha) = \sup\{\zeta:\alpha$ is $\zeta$-Mahlo$\}$.
If the set $\{i < \mu : \alpha_i \text{ is not $\alpha_i$-Mahlo}\}$ is
stationary then as before, for some $\gamma \in N$ we have
$\{i : F_{1,m(0)}(\alpha_i) = \gamma\}$ is stationary. Let
$F_{1,m(1)}(\alpha)$ --- a club of $\alpha$ such 
that if $\alpha$ is not $(\gamma+1)$-Mahlo then the club has no $\gamma$-Mahlo member.
Finish as in the proof of {Fact} $(\delta)$. 
\end{PROOF} 

\sn
Together we are done.
\end{PROOF}

\begin{remark}\label{e17} 
We can continue, and say more.
\end{remark}

\begin{lemma}\label{e20} 
$1)$ Suppose $\lambda > \mu > \theta$ are 
regular cardinals, $n\ge 2$, and
\begin{enumerate}
    \item[$(i)$] For every regular cardinal $\kappa$, if $\lambda > \kappa \ge \theta$ then $\kappa \not\to [\theta]^{<\omega}_{\sigma(1)}$.
\sn
    \item[$(ii)$] For some $\alpha(*) < \mu$, for every regular
    $\kappa \in (\alpha(*),\lambda)$, $\kappa \not\to [\alpha(*)]^n_{\sigma(2)}$.
\end{enumerate}
\underline{Then}
\begin{enumerate} 
    \item[$(a)$] $\lambda \not\to [\mu]^{n+1}_\sigma$, where 
    $\sigma = \min\{\sigma(1),\sigma(2)\}$.
\sn
    \item[$(b)$] There are functions $d_2 : [\lambda]^{n+1} \to \sigma(2)$ and
    $d_1 : [\lambda]^3 \to \sigma(1)$ such that for every $W \in [\lambda]^\mu$
    we have $d''_1([W]^3) = \sigma(1)$ or $d''_2([W]^{n+1}) = \sigma(2)$.
\end{enumerate}

\sn 
$2)$ Suppose $\lambda > \mu>\theta$ are regular cardinals, and
\begin{enumerate} 
    \item[$(i)$] For every regular $\kappa \in [\theta,\lambda)$ we have 
    $\kappa \not\to [\theta]^{<\omega}_{\sigma(1)}$. 
\sn
    \item[($ii$)] $\sup\{\kappa < \lambda : \kappa \text{ regular}\} \not\to [\mu]^n_{\sigma(2)}$.
\end{enumerate}
Then
\begin{enumerate} 
    \item[$(a)$] $\lambda\not\to [\mu]^{2n}_\sigma$, where 
    $\sigma = \min\{\sigma(1),\sigma(2)\}$.
\sn
    \item[$(b)$] There are functions $d_1 : [\lambda]^3 \to \sigma(1)$,
    $d_2 : [\lambda]^{2n} \to \sigma(2)$ such that for every 
    $W \in [\lambda]^\mu$ we have $d''_1([W]^3) = \sigma(1)$ or 
    $d''_2([W]^{2n}) = \sigma(2)$.
\end{enumerate}
\end{lemma}

The proof is similar to that of \cite[3.3,3.2]{Sh:276}.

\begin{PROOF}{\ref{e20}} 
1) For each $i$, $0 < i < \lambda_i$, we choose 
$C_i$ such that if $i$ is a successor ordinal then $C_i = \{i-1,0\}$,
and if $i$ is a limit ordinal then $C_i$ is a club of $i$ of order type
$\cf(i)$ containing 0 such that 
$[\cf(i) < i \Rightarrow \cf(i) < \min(C_i \setminus \{0\})$] 
and $C_i\setminus \acc(C_i)$\
contains only successor ordinals. 

Now for $\alpha<\beta$, $\alpha > 0$ we define $\gamma^+_\ell(\beta,\alpha)$, $\gamma^\supminus_\ell(\beta,\alpha)$ by induction on $\ell$, and then
$\kappa(\beta,\alpha),\,\eps(\beta,\alpha)$. 
\begin{enumerate}
    \item[(A)] $\gamma^+_0(\beta,\alpha) = \beta$, $\gamma^\supminus_0(\beta,\alpha ) = 0$.
\sn
    \item[(B)] If $\gamma^+_\ell(\beta ,\alpha )$ is defined and $> \alpha$ and
    $\alpha$ is not an accumulation point of $C_{\gamma^+_\ell(\beta ,\alpha )}$ then we let $\gamma^\supminus_{\ell+1}(\beta ,\alpha )$ be the maximal member of $C_{\gamma^+_\ell(\beta,\alpha)}$ which is $< \alpha$ and $\gamma^+_{\ell+1}(\beta,\alpha)$ is the minimal member of $C_{\gamma^+_\ell(\beta, \alpha)}$ which is $\ge \alpha$ 
    (by the choice of $C_{\gamma^+_\ell(\beta ,\alpha )}$ and the demands on $\gamma^+_\ell(\beta ,\alpha)$ they are well defined). 

\hspace{-1.6cm} So
\sn
    \item[(B1)] 
    \begin{enumerate}
        \item $\gamma^\supminus_\ell(\beta ,\alpha ) < \alpha \le \gamma^+_\ell(\beta,\alpha )$, and if the equality holds then $\gamma^+_{\ell+1}(\beta ,\alpha)$ is not defined.
\sn
        \item $\gamma^+_{\ell+1}(\beta ,\alpha) < \gamma^+_\ell(\beta ,\alpha)$ when both are defined. 
    \end{enumerate}
\sn
    \item[(C)] Let $k = k(\beta ,\alpha)$ be the maximal number $k$ such that $\gamma^+_k(\beta ,\alpha )$ is defined (it is well defined as 
    $\LL \gamma^+_\ell(\beta, \alpha) : \ell < \omega \RR$ is strictly decreasing). So
\sn
    \item[(C1)] $\gamma^+_{k(\beta,\alpha)}(\beta,\alpha) = \alpha$ or $\gamma^+_{k(\beta,\alpha)} > \alpha$, $\gamma^+_{k(\beta,\alpha)}$ is a limit ordinal and $\alpha$ is an accumulation point of $C_{\gamma^+_{k(\beta,\alpha)}}(\beta,\alpha)$. 
\sn
    \item[(D)] For $m \le k(\beta ,\alpha )$ let us define 
    $$\eps_m(\beta ,\alpha ) = \max\{\gamma^\supminus_\ell(\beta ,\alpha ) + 1 : \ell \le m\}.$$
    
\hspace{-1.6cm}  Note

    \item[(D1)]
    \begin{enumerate}
        \item $\eps_m(\beta ,\alpha ) \le \alpha$ (if defined).
\sn
        \item If $\alpha$ is limit then $\eps_m(\beta ,\alpha ) < \alpha$ (if defined).
\sn
        \item If $\eps_m(\beta ,\alpha) \le \xi \le \alpha$ \underline{then} for every $\ell \le m$ we have 
        $$\gamma^+_\ell(\beta ,\alpha ) = \gamma^+_\ell(\beta ,\xi ),\quad\gamma^\supminus_\ell(\beta ,\alpha ) = \gamma^\supminus_\ell(\beta ,\xi),\quad \eps_\ell(\beta ,\alpha ) = \eps_\ell (\beta ,\xi).$$ 
        (Explanation for (c): if $\eps_m(\beta ,\alpha ) < \alpha$ this is easy (check the definition) and if $\eps_m(\beta ,\alpha ) = \alpha$, necessarily $\xi = \alpha$ and it is trivial.)
\sn
        \item If $\ell \le m$ then $\eps_\ell (\beta ,\alpha) \le \eps_m(\beta,\alpha)$.
    \end{enumerate} 
\end{enumerate}
For a regular $\kappa \in (\alpha(*),\lambda)$ let 
$g^1_\kappa : [\kappa]^{<\omega} \to \sigma(2)$ exemplify 
$\kappa \not\to [\theta]^{<\omega}_{\sigma(1)}$, and
for every regular cardinal $\kappa \in [\theta,\lambda)$ let
$g^2_\kappa : [\kappa]^n \to \sigma(2)$ 
exemplify $\kappa \not\to [\alpha (*)]^n_{\sigma(2)}$. 

\sn
Let us define the colourings: 

Let $\alpha_0 > \alpha_1 > \ldots > \alpha_n$. (Remember $n\ge 2$.)

Let $n = n(\alpha_0,\alpha_1,\alpha_2)$ be the maximal natural number such that:
\begin{enumerate}
    \item[(i)] $\eps_n(\alpha_0,\alpha_1) < \alpha_0$ is well defined.
\sn
    \item[(ii)] $\gamma^\supminus_\ell(\alpha_0,\alpha_1) = \gamma^\supminus_\ell(\alpha_0,\alpha_2)$ for $\ell \le n$.
\end{enumerate}

We define $d_2(\alpha_0,\alpha_1,\ldots ,\alpha_n)$ as 
$g^2_\kappa(\beta_1,\ldots ,\beta_n)$, where 
$$
\begin{aligned}
\kappa =&\ \cf\big(\gamma^+_{n(\alpha_0,\alpha_1,\alpha_2)}(\alpha_0,\alpha_1)\big),\cr
\beta_\ell =&\ \otp \big(\alpha_\ell\cap
C_{\gamma^+_{n(\alpha_0,\alpha_1,\alpha_2)}(\alpha_0,\alpha_1)}\big).
\end{aligned}
$$ 

\sn
Next we define $d_1(\alpha_0,\alpha_1,\alpha_2)$ . 

Let $i(*) = \sup \big(C_{\gamma^+_n(\alpha_0,\alpha_2)} \cap
C_{\gamma^+_n(\alpha_1,\alpha_2)}\big)$, where 
$n = n(\alpha_0,\alpha_1,\alpha_2)$. Let $E$ be the equivalence 
relation on $C_{\gamma^+_n(\alpha_0,\alpha_1)} \setminus i(*)$
defined by 
$$\gamma_1\ E\ \gamma_2 \Leftrightarrow \big( \forall \gamma \in
C_{\gamma^+_n(\alpha_0,\alpha_2)} \big) [\gamma_1 < \gamma \Leftrightarrow \gamma_2 < \gamma ].$$ 

If the set $w = \big\{\gamma \in C_{\gamma^+_n(\alpha_0,\alpha_1)} :\gamma > i(*),\ \gamma = \min \gamma /E\big\}$ is finite, we let
$d_1(\alpha_0,\alpha_1,\alpha_2)$ be 
$g^1_\kappa \big(\{\beta_\gamma : \gamma \in w\}\big)$, where 
$\kappa = \big|C_{\gamma^+_n(\alpha_0,\alpha_1)}\big|$ and
$$\beta_\gamma = \otp\big(\gamma \cap 
C_{\gamma^+_n(\alpha_0,\alpha_1)}\big).$$

We have defined $d_1$, $d_2$ required in condition $(b)$
(though have not yet proved that they work) 
We still have to define $d$ (exemplifying $\lambda\not\to
[\mu]^{n+1}_\ell)$. Let $n \ge 3$: for $\alpha_0 > \alpha_1 >\ldots > \alpha_n$,
we let $d(\alpha_0,\ldots,\alpha_n)$ be $d_1(\alpha_0,\alpha_1,\alpha_2)$ if 
$w$ defined during the definition has odd number of members 
and $d_2(\alpha_0,\ldots
,\alpha_n)$ otherwise. 

Now suppose $Y$ is a subset of $\lambda$ of order type $\mu$, and
let $\delta = \sup Y$. Let $M$ be a model with universe $\lambda$ and with
relations $Y$ and $\{(i,j):i\in C_j\}$. Let $\LL N_i : i <
\mu \RR$ be an increasing continuous sequence of elementary
submodels of $M$ of cardinality $< \mu$ such that $\alpha(i) =
\alpha_i = \min(Y\setminus N_i)$ belongs to $N_{i+1}$, $\sup(N\cap \alpha_i) =
\sup(N\cap\delta)$. Let $N =\bigcup\limits_{i < \mu} N_i$. Let $\delta(i) =\delta_i =
\sup(N_i\cap \alpha_i)$, so $0 < \delta_i \le \alpha_i$, and let $n =
n_i$ be the first natural number such that $\delta_i$ an accumulation
point of $C^i = C_{\gamma^+_n(\alpha_i,\delta(i))}$, let $\eps_i =
\eps_{n(i)}(\alpha_i,\delta_i)$. Note that $\gamma^+_n(\alpha_i,\delta_i) = \gamma^+_n(\alpha_i,\eps_i)$ 
hence it belongs to $N$. 
 
\sn
\textbf{Case I}: { For some (limit) $i < \mu$, $\cf(i)
\ge\theta $ and $ (\forall \gamma < i)[\gamma + \alpha(*) < i]$ such that for
arbitrarily large $j < i$, $C^i\cap N_j$ is bounded in
$N_j\cap \delta=N_j\cap\delta_j$.}\\
\noindent 
This is just
like the last part in the proof of \cite[3.3]{Sh:276}, using $g^1_\kappa$ and $d_1$
for $\kappa=\cf(\gamma^+_{n_i}(\alpha_i,\delta_i))$. 

\sn
\textbf{Case II}: {Not case I}.

Let $S_0 = \{i < \mu : (\forall \alpha < i)[\gamma + \alpha(*) < i]$, 
$\cf(i) = \theta\}$. So for every $i \in S_0$, for some $j(i) < i$,
$$(\forall j) \big[j\in (j(i),i)\Rightarrow C^i\cap N_j \text{ is unbounded in }\delta_j\big].$$ 
But as $C^i \cap \delta_i$ is a club of $\delta_i$,
clearly $(\forall j)\big[j \in (j(i),i)\Rightarrow\delta_j\in C^i\big]$. 

We can also demand $j(i) > \eps_{n(\alpha(i),\delta(i))}(\alpha(i),\delta(i))$.

As $S_0$ is stationary, by `not {case I},' for some stationary
$S_1 \subseteq S_0$ and $n(*)$, $j(*)$ we have $(\forall i \in
S_1)\big[j(i) = j(*)\wedge n(\alpha(i),\delta_i) = n(*)\big]$. 

Choose $i(*)\in S_1$, $i(*) = \sup(i(*)\cap S_1)$, such that the order
type of $S_1\cap i(*)$ is $i(*)> \alpha(*)$. Now if $i_2 < i_1\in
S_1\cap i(*)$ then $n(\alpha_{i(*)},\alpha_{i_1},\alpha_{i_2}) = n(*)$.
Now $L_{i(*)} = \big\{\otp(\alpha_i \cap C^{i(*)}):i \in
S_1\cap i(*)\big\}$ are pairwise distinct and are ordinals $<$
$\kappa = |C^{i(*)}|$, and the set has order type $\alpha(*)$.
Now apply the definitions of $d_2$ and $g^2_\kappa$ on $L_{i(*)}$. 
\medskip
\noindent 2) The proof is like the proof of part (1),but for 
$\alpha_0 > \alpha_1 >\cdots$ we let
$d_2(\alpha_0,\ldots ,\alpha_{2n-1}) = g^2_\kappa(\beta_0,\ldots,\beta_n)$, 
where 
$$\beta_\ell =
\otp\big(C_{\gamma^+_n(\beta_{2\ell},\beta_{2\ell+1})}(\beta_{2\ell},\beta_{2\ell+1})
\cap\beta_{2\ell+1}\big)$$
and in case II note that the analysis gives $\mu$
possible $\beta_\ell$-s so that we can apply the definition of
$g^2_\kappa$. 
\end{PROOF}

\begin{definition}\label{e23} 
Let $\lambda \not\to_\stg [\mu]^n_\theta$ mean: if 
$d : [\lambda]^n \to \theta$, $\LL \alpha_i : i < \mu \RR$ is
strictly increasing continuous, and for $i < j < \mu$, 
$\gamma_{i,j} \in [\alpha_i,\alpha_{i+1})$ then 
$$\theta = \big\{d(w) : \hbox{for some } j < \mu,\,\, w \in
\big[\{\gamma_{i,j} : i < j\}\big]^n\big\}.$$
\end{definition}

\begin{lemma}\label{e26} 
1) $\aleph_t \not\to [\aleph_1]^{n+1}_{\aleph_0}$ for $n \ge 1$.

\sn
2) $\aleph_n \not\to_\stg [\aleph_1]^{n+1}_{\aleph_0}$ for $n \ge 1$.
\end{lemma}

\begin{PROOF}{\ref{e26}} 
1) For $n=2$ this is a theorem of Todor{\v c}evi{\v c} \cite{To2},
and if it holds for $n\ge 2$ by \ref{e20}(1) 
we get that it holds for n+1 (with $n$, $\lambda$, $\mu$, $\theta$, $\alpha(*)$, $\sigma(1)$, $\sigma(2)$ there corresponding to 
$n+1$, $\aleph_{n+1}$, $\aleph_1$, $\aleph_0$, $\aleph_0$,
$\aleph_0$, $\aleph_0$ here). 

\sn
2) Similar. 
\end{PROOF}

\bibliographystyle{amsalpha}
\bibliography{shlhetal}
\end{document}